\documentclass[mnsc,nonblindrev]{formatfile} % current default for manuscript submission

\OneAndAHalfSpacedXI

\usepackage[all]{xypic}
\usepackage{graphics}
\usepackage{epstopdf}
\usepackage{color}
\usepackage{natbib}
 \bibpunct[, ]{(}{)}{,}{a}{}{,}%
 \def\bibfont{\small}%
 \def\bibsep{\smallskipamount}%
 \def\bibhang{24pt}%
\TheoremsNumberedThrough     % Preferred (Theorem 1, Lemma 1, Theorem 2)
\ECRepeatTheorems

\EquationsNumberedThrough    % Default: (1), (2), ...
\MANUSCRIPTNO{-}

\begin{document}
%%%%%%%%%%%%%%%%

 \RUNAUTHOR{Faghih, Roozbehani, Dahleh} % for four or more authors

\RUNTITLE{On the Economic Value and Price-Responsiveness of Ramp-Constrained Storage}

% \TITLE{Bundling Information Goods of Decreasing Value}
\TITLE{On the Economic Value and Price-Responsiveness of Ramp-Constrained Storage}

\ARTICLEAUTHORS{%
\AUTHOR{Ali Faghih}
\AFF{ (Corresponding Author, Email: \EMAIL{afaghih@mit.edu})\\Massachusetts Institute of Technology (MIT), 77 Mass Ave, Cambridge, MA, USA 02139 \\ Electrical Engineering and Computer Science (EECS) Department, Laboratory for Information and Decision Systems (LIDS)}  %, \URL{}}
\AUTHOR{Mardavij Roozbehani}
\AFF{LIDS, MIT, Cambridge, MA, USA 02139. \EMAIL{mardavij@mit.edu}}
\AUTHOR{Munther A. Dahleh}
\AFF{EECS, LIDS, MIT, Cambridge, MA, USA 02139. \EMAIL{dahleh@mit.edu}}
} % end of the block

\ABSTRACT{%
The primary concerns of this paper are twofold: to understand the economic value of storage in the presence of ramp constraints and exogenous electricity prices, and to understand the implications of the associated optimal storage management policy on qualitative and quantitative characteristics of storage response to real-time prices. We present an analytic characterization of the optimal policy, along with the associated finite-horizon time-averaged value of storage. We also derive an analytical upperbound on the infinite-horizon time-averaged value of storage. This bound is valid for any achievable realization of prices when the support of the distribution is fixed, and highlights the dependence of the value of storage on ramp constraints and storage capacity. While the value of storage is a non-decreasing function of price volatility, due to the finite ramp rate, the value of storage saturates quickly as the capacity increases, regardless of volatility. To study the implications of the optimal policy, we first present computational experiments that suggest that optimal utilization of storage can, in expectation, induce a considerable amount of price elasticity near the average price, but little or no elasticity far from it. We then present a computational framework for understanding the behavior of storage as a function of price and the amount of stored energy, and for characterization of the buy/sell phase transition region in the price-state plane. Finally, we study the impact of market-based operation of storage on the required reserves, and show that the reserves may need to be expanded to accommodate market-based storage. 
\vspace{5pt}

%\noindent \emph{Funding Sources:} National Science Foundation Graduate Research Fellowship; MIT Energy Initiative Seed Fund; Draper Laboratories. The funding sources were not involved in conducting this research, preparing the article, or the decision to submit the article for publication

% Enter your abstract
}%




\maketitle
%%%%%%%%%%%%%%%%%%%%%%%%%%%%%%%%%%%%%%%%%%%%%%%%%%%%%%%%%%%%%%%%%%%%%%%%%%%%%%%%%%%%%%%%%%%%%%%%%%%%%%%%%%%%%%%%%%%%%%%%%%%%

\vspace{-15pt}
\section{Introduction} \label{introd}
The growing demand for electricity and the urge to reduce greenhouse emissions promote large-scale integration of renewable energy sources, as well as new storage and demand-response technologies to improve energy efficiency in the future grid. However, renewable energy sources are highly uncertain and intermittent. While energy storage technologies can help mitigate the intermittency and narrow the gap between generation from renewable resources and consumption, they may also add to the uncertainty in the system since the optimal response of storage to prices is a complicated function of both price and the amount of stored energy. Modeling and understanding the behavior of storage in response to real-time market prices is therefore critical for reliable operation of power systems with large amounts of storage.

Energy storage has a clear environmental value; it helps mitigate the intermittency of the renewable resources and thereby maximize their utilization. With sufficient storage capacity, there is no need to curtail generation from renewable energy sources when there is excess of it. Moreover, with proper control policies, storage can help incorporate more renewable resources without compromising reliability. Access to storage also reduces the risk associated with making advance commitments for a renewable generation owner. For instance, if the energy generated by the renewable source falls short of the committed level, the operator can compensate for the shortfall by extracting from the storage. Despite all these potential advantages of storage, if the economic value of storage as an arbitrage mechanism is not attractive, the markets may not invest sufficiently in storage. Hence, unless proper incentives and pricing policies are in place, the environmental and reliability values of storage might not materialize due to underinvestment. Therefore, there is a need for development of econometric models and characterization of the associated optimal policies that can be used for assessing the economic value of storage. This paper seeks to provide such characterization by presenting a model for optimal utilization of ramp-constrained storage in response to stochastically varying electricity prices.

Availability of econometric models of storage and characterization of the effects of storage on the overall price elasticity of demand is also important for system operators who need to maintain stability, and guarantee reliability of the system. It was shown in \cite{RoozbehaniTPS2011} that in power grids with information asymmetry between consumers, producers, and system operators, robustness of the system to disturbances is greatly affected by consumers' real-time valuation of electricity, and their response to real-time prices. It was shown that under real-time pricing, high price volatility can be associated with uncertainty in demand response and high price elasticity of demand (PED).

The existing literature covering various dimensions of storage management is extensive. The two major streams of literature on storage that are closely related to this work are the commodity trading/warehousing literature and the electricity storage literature. 
%%%%%commodity storage literature
The warehouse problem is a classical problem in the trading and commercial management of commodities, and has been studied in the literature extensively. As early as 1948, \cite{Cahn} introduced the problem of optimizing purchase, i.e. injecting into storage, and sale, i.e., withdrawing from storage, for the case of a warehouse with fixed size and an initial stock of a certain commodity. \cite{Bellman} formulated the warehouse problem in a dynamic programming framework, and \cite{Drey} showed that for deterministic prices, the optimal policy at each stage is to either fill up the storage, or empty it, or do nothing, depending on the stage price. \cite{Charnes} solved the warehouse problem for the case of stochastic prices and showed that the optimal policy is the same in the stochastic case (either fill up the storage, or empty it, or do nothing). However, in the warehouse problem studied by these works, there is no limit on the amount that can be injected into, or withdrawn from storage at each stage. 

\cite{Remp}, and more recently, \cite{Secomandi}, extended these results for the discrete-time case by imposing a limit on the amount that can be injected into or withdrawn from storage. \cite{Kam} solved the same extension of the warehouse problem in continuous time. Some other works have also considered this extension of the warehouse problem, including \cite{Wu}, who propose a heuristic for optimization of seasonal energy storage operations in the presence of ramp constraints, as well as \cite{Deval} who seek to optimize procurement, processing, and trading of commodities in a multi-period setting. Although the above-mentioned references and our work share similarities in the structure of the optimal policy and the associated value function, the differences in the assumptions on the stochastic price process make the analytical results of these papers different from one another. Unlike all the previous works in this area, we assume that the price at each time period is independent of previous prices. Under this assumption we derive explicit formulas (recursive and/or closed-form) for the thresholds of the optimal policy and for the economic value of storage. We justify our assumption on prices by testing the performance of our optimal policy against price data from real-time markets of the \cite{PJM} and \cite{isone}. Another justification for this assumption is that in practice, empirical estimation of conditional distributions (needed for a Markovian price model) requires significant amounts of data. Back of the envelop calculations show that collecting this much data would require going too far back in the price history. However, due to non-stationarity, doing so would make the data irrelevant. Although one can resort to calibrated models for estimation of correlation in data, or try to learn the thresholds directly, we do not pursue these directions in this paper. In addition, we derive an upperbound on the optimal average value per stage of storage in the infinite-horizon case, explicitly highlighting its dependence on storage capacity and the ramp rate. Also different from previous works, we present a computational framework for estimating the average PED obtained from averaging out the stage-dependence and internal state of storage. We extend our study of the PED induced by storage by addressing the PED as a function of the internal state of storage. This is particularly important in the context of electricity markets because the aggregate price elasticity of demand can affect price volatility and sensitivity to disturbances (\citealt{RoozbehaniTPS2011}). We also highlight another important aspect of state dependence of storage response by showing that under certain assumptions, if the system operator does not know the exact state of storage, the reserves may need to be expanded to accommodate market-based operation of storage.

%%%%%The electricity storage literature
The literature focusing particularly on electricity storage is also extensive.
\cite{Bannister} studied optimal management of a single storage connected to a general linear memoryless system in the presence of ramp constraints. However, in their model, the objective function is deterministic and the cost is known a priori. Also, \cite{Lee2} considered industrial consumers with time-of-use rates and used dynamic programming to determine optimal contracts and optimal sizes of battery storage systems for such consumers. Their work, like ours, pays special attention to the economic value of storage; however, they use a deterministic approach and relax ramp constraints. Several other works have studied the impacts of energy storage on the economics of integration of renewable sources. A renewable generation owner can connect a storage device to the renewable source to optimize the overall profit over time by deciding how much energy to commit to sell at each stage in the time horizon. Some more recent works such as \cite{Brown}, \cite{Gonz}, and \cite{Korp} approach this problem by deterministically solving this optimization problem for particular finite sample paths and then averaging the results of these paths. Their approach, as mentioned in \cite{Harsha} and \cite{Kim}, does not give an optimal policy that can be used in practice because their policy depends on the sample path. 

On the other hand, some recent works such as \cite{Bitar}, \cite{Harsha}, and \cite{Kim} use dynamic programming to address the problem of managing the revenue of a renewable generator using storage. In particular, \cite{Harsha} use a stochastic dynamic programming approach to study the optimal storage investment problem through characterization of optimal sizing of energy storage for efficient integration of renewable resources. In contrast to \cite{Harsha}, we explicitly include ramp constraints in our model, and highlight the effects of ramp constraints on the value of storage. Another contribution in this line of research is due to \cite{Kim}, who study the problem of making advance commitments for a wind generator in the presence of storage capacity and conversion losses. \cite{Kim} use dynamic programming to obtain the optimal commitment policy when the storage device can have conversion losses but under the assumption of uniformly distributed generation from the wind farm. This assumption allows them to derive the stationary distribution of the storage level
and use it to characterize the economic value of storage. In contrast, our results characterize the value of storage purely as an arbitrage mechanism that interacts only with the main grid with a guaranteed supply, without any particular assumption on price distribution other than the assumption of independent prices. 

The storage problem in this paper also has similarities with the inventory management literature in terms of the underlying dynamic programming problem and the corresponding optimal policy (for instance, see \citealt{Feder}, \citealt{Kap}, and \citealt{Goel}). However, the main difference comes from the fact that the inventory literature focus on optimally managing inventory in the presence of demand, while the model in this paper assumes that the main grid buys all the energy that we decide to supply, and instead, focuses on trading inventory under ramp constraints and maximizing profit by taking advantage of the fluctuations in the spot prices. Our setup also has similarities with the literature on reservoir hydroelectric system management (see, for instance, \citealt{Drouin} and \citealt{Lam}) especially those with price uncertainty. However, since the focus of our paper is on the economic value and price-responsiveness of electricity storage as an arbitrage mechanism under ramp constraints, the formulation of our model and our results are different from this line of research. 

The model set forth in this paper and some preliminary results were reported in \cite{Faghih}. This paper provides a comprehensive exposition which adds several new ideas, core analytical results, and systematic computational experiments. The contributions of this paper are summarized as follows:
First, we propose a dynamic model for optimal utilization of storage in the presence of ramp constraints under the assumption of independent and exogenous prices. This model assumes limited storage capacity and allows sell-back of energy to the grid. Using the principles of stochastic dynamic programming, we analytically characterize the optimal policy and the corresponding value function for the finite-horizon case. In particular, we provide recursive equations for computation of the exact value of storage for the finite-horizon storage problem. To verify the validity of our assumptions, we apply our finite-horizon optimal policy to real-time price data taken from the PJM Interconnection and the Independent System Operator (ISO) of New England. We define the \emph{competitive ratio} (CR) as the ratio of the value obtained from our optimal policy to the absolute maximum value that would have been obtained deterministically, had we known the entire price process a priori. We then show that the value obtained from our policy under the assumption of independent prices yields a relatively high CR when applied to real-world price data, reaching a CR of about $90\%$ in some cases. This suggests that the sensitivity of the optimal policy and value of storage to the assumption of independent prices is low. We then obtain a closed-form upperbound on the infinite-horizon optimal average value per stage of storage over all possible realizations of prices within a bounded support. To the best of our knowledge, this paper is the first to derive an analytical upperbound on the long-term average value of ramp-constrained storage. This result highlights how the capacity limit and the ramp constraint bound the value of storage. Next, we show that while the economic value of storage is a non-decreasing function of price volatility, the value of storage saturates quickly due to finite ramping rates as the capacity increases, regardless of price volatility. Our results on the economic value of storage can be used to evaluate certain decisions about investment in storage. 

Next, in order to study the average price elasticity of a storage system in an electricity market, we study the average PED in a simulated electricity market over a one-day time-horizon, where the term ``average PED" reflects the fact that the dependence of storage response on the stage and the internal state of storage has been averaged out. We show that optimal utilization of storage may, in expectation, induce a considerable amount of price elasticity near the average price, but little or no elasticity elsewhere. While the demand for electricity has often been considered to be highly inelastic, the existing literature on price elasticity are mostly based on empirical evidence and qualitative reasoning, see, for instance, \cite{Kirschen}, \cite{Kirschen2}, \cite{Yusta}, and \cite{Faruqui}. In this paper, we study price elasticity in a quantitative framework. To the best of our knowledge, this paper is the first to characterize the PED induced by storage through an input-output model of response to prices based on optimal control policies in the presence of ramp constraints. Finally, to examine how the storage response would have been had we not averaged out the state, we characterize the response of storage to exogenous prices and its dependence on the internal state of storage, and highlight the interplay between state-dependence and price-dependence of the response in a computational framework. To eliminate time-dependence, we study these relations in an infinite-horizon setting. We use policy iteration to characterize price responsiveness of a storage system over an infinite time-horizon as a function of the storage state, and characterize the buy/sell phase transition region in the price-state plane. To highlight an interesting implication of the price-state interplay, we study the impact of market-based operation of storage on the required reserves, and show that if the ISO does not have perfect information about the exact value of the storage state, the reserves may need to be expanded to accommodate market-based operation of storage. 

The remainder of this paper is organized as follows: In Section \ref{sec:model}, we introduce the dynamic model of utilization of storage. In Section \ref{sec:results}, we present the optimal policies for the storage management problem and the corresponding value function, and give an evaluation of the performance of the optimal policy. Next, in Section \ref{valchar} we derive an analytical upperbound on the optimal average value per stage of storage, and report our computational findings on the economic value of storage. We then discuss the implications of the optimal policy in Section \ref{sec:policyimp}, by studying the average PED of a storage system in Subsection \ref{sec:elasticity}, the price responsiveness of a storage system as a function of the storage state in Subsection \ref{infelast}, and the impact of market-based operation of storage on the required reserves in Subsection \ref{sec:reserves}. We conclude in Section \ref{sec:conclusions}.
%\vspace{-1pt}
%%%%%%%%%%%%%%%%%%%%%%%%%%%%%%%%%%%%%%%%%%%%%%%%%%%%%%%%%%%%%%%%%%%%%%%%%%%
%%%%%%%%%%%%%%%%%%%%%%%%%%%%%%%%%%%%%%%%%%%%%%%%%%%%%%%%%%%%%%%%%%%%%%%%%%%
\section{A Dynamic Model of Storage} \label{sec:model}

\subsection{Notation}
%%%%%%%%%%%%%%%%%%%%%%%%%%%%%%%%%%%%%%%%%%%%%%%%%%%%%%%%%%%%%%%%%%%%%%%%%%%
%%%%%%%%%%%%%%%%%%%%%%%%%%%%%%%%%%%%%%%%%%%%%%%%%%%%%%%%%%%%%%%%%%%%%%%%%%%
The set of positive real numbers (integers) is denoted by $\mathbb{R}_{+}$
($\mathbb{Z}_{+}$)$,$ and non-negative real numbers (integers) by $\overline
{\mathbb{R}}_{+}$ ($\overline{\mathbb{Z}}_{+}$). The probability mass function (PMF) of a random variable $\Lambda$ is denoted by $P_\Lambda$, and the cumulative distribution function (CDF) is denoted by $F_\Lambda$. We will simply use $P$ and $F$ when there is no ambiguity.%\vspace{-3pt}
%%%%%%%%%%%%%%%%%%%%%%%%%%%%%%%%%%%%%%%%%%%%%%%%%%%%%%%%%%%%%%%%%%%%%%%%%%%
%%%%%%%%%%%%%%%%%%%%%%%%%%%%%%%%%%%%%%%%%%%%%%%%%%%%%%%%%%%%%%%%%%%%%%%%%%%
\subsection{The Model} \label{sec:storFormu}

In this section, we develop a dynamic model for optimal management of ramp-constrained storage in the presence of stochastically-varying prices. We start by formulating the storage management problem as a stochastic dynamic programming problem over a finite horizon. 

%%%%%%%%%%%%%%%%%%%%%%%%%%%%%%%%%%%%%%%%%%%%%%%%%%%%%%%%%%%%%%%%%%%%%%%%%%%
%%%%%%%%%%%%%%%%%%%%%%%%%%%%%%%%%%%%%%%%%%%%%%%%%%%%%%%%%%%%%%%%%%%%%%%%%%%
\subsubsection*{The Decisions:} The decision set of the storage owner or consumer at each discrete instant of time $k\in\overline{\mathbb{Z}}_{+}$ is characterized by a pair
\vspace{-1pt}
\begin{equation}
(v_{k}^{\text{in}},v_{k}^{\text{out}})\in \lbrack0,\overline{v}^{\text{in}}]\times\lbrack0,\overline{v}^{\text{out}}]
\label{decisionsconstraints}%
\end{equation}
%\vspace{-10pt}
\noindent where, $v_{k}^{\text{in}}$ and $v_{k}^{\text{out}}$ are, respectively, the amount of power that the consumer injects in, or withdraws from the storage. The corresponding upper bounds ($\overline{v}^{\text{in}}$ and $\overline{v}^{\text{out}}$) represent the physical ramp constraints on storage. Also, $v_{k}=v_{k}^{\text{in}}-v_{k}^{\text{out}} \in \lbrack-\overline{v}^{\text{out}},\overline{v}^{\text{in}}]$ denotes the net storage response. With a slight abuse of terminology, we may also refer to the storage response $v_k$ as demand, with the understanding that $v_k \leq 0$ implies a negative demand.
\subsubsection*{The Price:} The price at each stage ${\mathbf \lambda_k}$ is sampled from an exogenous stochastic process that is independently distributed across time, with mean $\overline\lambda_{k}$ and standard deviation $\sigma_k$, and support over $[\lambda^{\min}_{k},\lambda^{\max}_{k}]\subset [0,\infty)$. It is assumed that at the beginning of each time interval $\left[  k,k+1\right] ,$ the random variable $\lambda_{k}\,$ is materialized and revealed to the consumer. Note that the distributions of $\lambda_{k}$ can be different at different stages; however, we assume that the price distribution at each stage is known a priori. We also assume that the feed-in and usage tariffs are the same, i.e., $\lambda_{k}$ is the price per unit for
both purchase (corresponding to $v_{k}\geq0$)\ and sell-back
(corresponding to $v_{k}\leq0$), and there are no transaction costs.
\subsubsection*{The States:}
The storage state is characterized by a variable
\vspace{-8pt}
\begin{equation}
s_{k}  \in \left[  0,\overline{s}\right]
\label{thestates}%
\end{equation}
%\vspace{-8pt}
\noindent where $s_{k}$ is the amount of energy stored, and $\overline
{s}$ is the upper bound on storage capacity. The state $s_{k}$ evolves according to:
\vspace{-8pt}
\begin{align}
s_{k+1}  &  =\beta s_{k}+\eta^{\text{in}}v_{k}^{\text{in}}-\eta^{\text{out}
}v_{k}^{\text{out}} \label{storageevolution}
\end{align}
%\vspace{-5pt}
\noindent where $\beta\leq1$ is the decay factor, $\eta^{\text{in}}\leq1$
and $\eta^{\text{out}}\geq1$ are charging and discharging efficiency
factors. Note that efficiency factors and ramp rates might in general be
complicated functions of the operating point, i.e., the storage level, but in this work we focus our attention on an ideal case. The
idealized model of the dynamics of storage can be written as:
\vspace{-3pt}
\begin{equation}
s_{k+1}=s_{k}+v_{k},\qquad v_{k}\in\lbrack-\overline{v}_{\text{out}}
,\overline{v}_{\text{in}}] \label{storageevolution2}
\vspace{-5pt}
\end{equation}
which corresponds to $\beta=1,$ $\eta^{\text{in}}=1,$ and $\eta^{\text{out}
}=1.$
%\vspace{-2pt}
\subsubsection*{Penalty and Salvage Value:} There is a penalty $h_{k}(s_{k})$ associated with storage, where the sequence of functions $h_{k}:\overline{\mathbb{R}}_{+}\mapsto\overline{\mathbb{R}}_{+}$ are assumed to be monotonic. Also, a salvage value of $\hat{t} \in [\lambda^{\min}_{N},\lambda^{\max}_{N}]$ is assigned to each unit of energy left in storage by the end of the time-horizon.
%\vspace{2pt}
\subsubsection*{The Optimization-Based Model of Ideal Storage:} Since our goal in this paper is to develop tractable models (and derive bounds) that effectively highlight the important structural features of the optimal control law and the associated economic value of storage with an emphasis on the ramp constraint and storage capacity, we will adopt the idealized model of storage. Nevertheless, in terms of methodology, it would be straightforward to include a ``price-adjustment" factor in our formulation to account for injection-withdrawal losses in a similar manner done in \cite{Secomandi}. In addition, the piecewise-linear penalty function that we have embedded into our model can be used as a surrogate cost to model the dissipation losses associated with keeping energy in storage. Note also that according to \cite{Qian}, the academic state of the art is around $95\%$ efficiency for a battery pack and $93\%$ for the overall system with converters. The industrial state of the art is around $90\%$ efficiency for batteries (\citealt{A123}). We also assume that the ramp constraint is symmetric, i.e. $\overline{v}^{\text{in}}=\overline{v}^{\text{out}}=\overline{v}$. The idealized storage management problem can be formulated as a finite-horizon dynamic programming problem as follows:
\begin{align}
&
\begin{array}
[c]{cl}
\min & \mathbf{E}\left[
{\displaystyle\sum\nolimits_{k=0}^{N-1}}
\{h_{k}(s_{k})+\lambda_{k}v_{k}\} -\hat{t}s_N\right] \vspace*{0.09in}
\end{array}
\label{storage}\\
&
\begin{array}
[c]{cl}
\text{s.t.\hspace{0.09in}} & s_{k+1}=s_{k}+v_{k}\vspace*{0.08in}\\
& s_{k}\in\lbrack0,\overline{s}]\vspace*{0.08in}\\
& v_{k}\in\lbrack-\overline{v},\overline{v}] \vspace*{0.08in}\\
& \lambda_{k}~ \text{exogenous, and independently}\\
& ~~~~\text{distributed according to a PMF }P_{k}   
\end{array}
\nonumber
\end{align}
\begin{remark}
We first formulate and solve the storage problem for the finite-horizon case. Later in Section \ref{analyticbound}, we consider the infinite-horizon case and obtain an upperbound on the optimal average value per stage of storage.%, 
\end{remark}
%\vspace{-3pt}
%%%%%%%%%%%%%%%%%%%%%%%%%%%%%%%%%%%%%%%%%%%%%%%%%%%%%%%%%%%%%%%%%%%%%%%%%%%%%%%%
\section{Characterization of the Optimal Policy and an Evaluation of Its Performance} \label{sec:results}
\subsection{The Optimal Policy}
In this subsection, we characterize the optimal policy and the value function for problem (\ref{storage}) based on principles of stochastic dynamic programming. 
\smallbreak
\noindent {\bf Definition 1.} Given a sequence of probability mass functions $P_{k},~k=0,1,\dots,N$, let $\Theta_{k}$ and $\psi_{k}$ be sequences of maps
from the set of all subsets of $\overline{\mathbb{R}}_{+}$ to $\overline
{\mathbb{R}}_{+}$, defined as follows:
\begin{align}
\Theta_{k} & :I\mapsto\sum\limits_{\theta\in I}\theta P_{k}\left(  \theta\right)
,\text{\qquad}\forall I\subset\overline{\mathbb{R}}_{+} \label{map1}\\
\psi_{k} & :I\mapsto\sum\limits_{\theta\in I}P_{k}\left(  \theta\right)
,\text{\qquad}\forall I\subset\overline{\mathbb{R}}_{+}.\label{map2}
\end{align}

\noindent Given $\overline{v}\in\mathbb{R}_{+},$ and maps $\Theta_{k}$ and $\psi_{k}$ as defined
in (\ref{map1}) and (\ref{map2}), $\Phi_{k}^{\overline{v}}$ is the map from the set of all subsets of
$\overline{\mathbb{R}}_{+}$ to $\mathbb{R}$ defined according to%
\vspace{-3pt}
\begin{equation}
\Phi_{k}^{\overline{v}}:I\mapsto\overline{v}\left(  \Theta_{k}-\rho\psi_{k}\right)  I,\text{\qquad
}\forall I\subset\overline{\mathbb{R}}_{+}, \label{phimap}
\vspace{-5pt}
\end{equation}
where
\vspace{-5pt}
\[
\rho=\inf I.
\vspace{-3pt}
\]
For instance, $\Phi_{k}^{1}$ maps an interval $\left(  a,b\right)  $ to $\left(
\Theta_{k}-a\psi_{k}\right)  \left(  a,b\right) .$

\begin{theorem}
\label{storageTHM}Consider the finite-horizon storage management problem (\ref{storage})
with $\overline{s}=n\overline{v}$ for some $n\in \mathbb{Z}_{+}$. Furthermore, assume that the penalty functions $h_{k}:[  0,\infty)  \rightarrow
\lbrack0,\infty),$ $k=0,\dots,N$ are piecewise linear non-decreasing convex functions of the
form:
\begin{equation}
h_{k}\left(  s\right)  =h^{i}_{k}s+c^{i}_{k},\text{ }s\in\left[  i\overline
{v},(i+1)\overline{v}\right),~i \in \overline{\mathbb{Z}}_{+}   \label{penal}
\vspace{-11pt}
\end{equation}
\noindent Then, \vspace{3pt}\newline \textrm{(i)} the optimal policy is characterized
as follows:\\
\vspace{2pt}
\noindent if $s_{k}\in\left[i\overline
{v},(i+1)\overline{v}\right)$, for some $i\in\overline{\mathbb{Z}}_{+}$, then %\{0,1,2,...,n-1\}
\begin{equation}
v_{k}^{\ast}=\left\{
\begin{array}
[c]{lclcl}%
%-s_{k} & , & t^{0}_{k+1}<\lambda_{k}&\text{for}&i=0\\[0.05in]
\max(-s_{k} ,-\overline{v}) & ,& t^{\max(0,i-1)}_{k+1}<\lambda_{k}&\text{for}&\text{all }i\\[0.05in]
i\overline{v}-s_{k} &,& t^{i}_{k+1}<\lambda_{k}\leq t^{i-1}_{k+1}&\text{for}&i\geq1\\[0.05in]
(i+1)\overline{v}-s_{k} &, & t^{i+1}_{k+1}<\lambda_{k}\leq t^{i}_{k+1}&\text{for}&\text{all }i\\[0.05in]
\overline{v} & , & \hspace*{0.45in}\lambda_{k}\leq t^{i+1}_{k+1}&\text{for}&\text{all }i
\end{array}
\right.
\vspace{3pt} \label{optpol}
\end{equation}
\vspace{8pt}
where the thresholds are computed via the following recursive equations:
\begin{equation}
 \vspace{10pt}
\begin{array}
[c]{lcl}
t^{i}_{N}=\hat{t},~~~~~~~~~i\in\{0,1,2,...,n-1\}\vspace{3pt}\\[0.05in]
t^{i}_{N}=-h^{i}_{N},~~~~~i\geq n  \vspace{12pt}\\[0.05in]
 \vspace{4pt}
\text{for}~k<N:\vspace{3pt}\\
t^{0}_{k}= t^{1}_{k+1} 
 + \Phi_{k}^{1}(t^{1}_{k+1}, \lambda^{\max}_{k}] -h^{0}_{k}\\[0.15in] %\!

%\!
t^{i}_{k}=t^{i-1}_{k+1}-h^{i}_{k}  +
 
 \Phi_{k}^{1}(t^{i+1}_{k+1}, t^{i-1}_{k+1}]+(t^{i+1}_{k+1}-t^{i-1}_{k+1})F_{k}\left(  t^{i-1}_{k+1}\right),~~~~~~~~~~i\geq 1  \\ 
\end{array}
\vspace{-6pt}
\label{thresh}
\end{equation}
\newline\textrm{(ii) }the value function is a piecewise linear convex function of the form:
%\vspace{-3pt}
\begin{equation}
V_{k}\left(  s\right)  =-t^{i}_{k}s+e^{i}_{k},\text{ }s\in\left[  i\overline
{v},(i+1)\overline{v}\right), i\in\overline{\mathbb{Z}}_{+} \label{valuefunc}
%\vspace{-3pt}
\end{equation}
where $t^{i}_{k}$ are the thresholds given in \eqref{thresh} and the intercepts $e^{i}_{k}$ are computed via the following recursive equations:
\begin{equation}
\vspace{5pt}
\begin{array}
[l]{lll}
\!e^{i}_{N}=0,~~~~~~~~~~~~~~i\in\{0,1,2,...,n-1\}\vspace{3pt}\\[0.06in]
\!e^{i}_{N}=\overline{s}(t^{i}_{N}-\hat{t}),\hspace{1pt}~~~~~i\geq n \vspace{-3pt}\\[0.06in]
\\
\text{for}~k<N:\vspace{5pt}\\
\!e^{0}_{k}\!=\!c^{0}_{k}+ e^{0}_{k+1}+(\overline{v}\lambda^{\min}_{k}+e^{1}_{k+1}-e^{0}_{k+1}
-\overline{v}t^{1}_{k+1})F_{k}(  t^{0}_{k+1})
  \vspace{12pt}\\
\!e_{k}^{i}\!=\! c_{k}^{i}\!+\!\overline{v}\overline{\lambda}+f(t^{i-1}_{k+1},t^i_{k+1},t^{i+1}_{k+1},e^{i-1}_{k+1},e^{i}_{k+1},e^{i+1}_{k+1})+\!g(t^{i-1}_{k+1},t^{i}_{k+1},t^{i+1}_{k+1}),~i\geq 1
\end{array}
\label{e_k}
\vspace{6pt}
\end{equation}
where the functions $f$ and $g$ are given by\vspace{3pt}
\[
\vspace{3pt}
\begin{array}
[l]{lll}
\hspace{-.65in}f(\cdot)= e^{i-1}_{k+1}-\overline{v}t_{k+1}^{i+1}+ (e^{i}_{k+1}-e^{i-1}_{k+1})F_{k}(t^{i-1}_{k+1})\vspace{0.0505in}+(e^{i+1}_{k+1}-e^{i}_{k+1})F_{k}(t^{i}_{k+1}) \vspace{0.15in}\\
\hspace{-.65in}g(\cdot)=\left(  i+1\right)\!  \Phi_{k}^{\overline{v}}(t_{k+1}^{i+1},t_{k+1}^{i}] +i\Phi_{k}^{\overline{v}}(t_{k+1}^{i},t_{k+1}^{i-1}]
 \vspace{0.05in}-\Phi_{k}^{\overline{v}}(t_{k+1}^{i-1},\lambda^{\max}_{k}] -\Phi_{k}^{\overline{v}%
}(t_{k+1}^{i+1},\lambda^{\max}_{k}].
\end{array}
\label{f_g}
\]

\end{theorem}

\proof{Proof.}
Please see the Appendix.
\Halmos
\endproof 
Figure \ref{threshpic} shows how the thresholds vary with time and state for the case of a discretized truncated log-normal distribution with mean $\overline{\lambda}_{k}=49$ and $\sigma_{k}=9$ for all $k$, i.e. for independently and identically distributed (i.i.d) prices. For generating this plot, we set $\overline{v}=1$, $n=15$, $N=24$, $\hat{t}=\overline{\lambda}_{k}$, and assume no storage penalties (i.e. $h^{i}_{k}=0$ for $i<n$ and all $k$). 

\begin{figure}
[ptbh]
\begin{center}
\vspace{-13pt}
\includegraphics[scale=0.5]{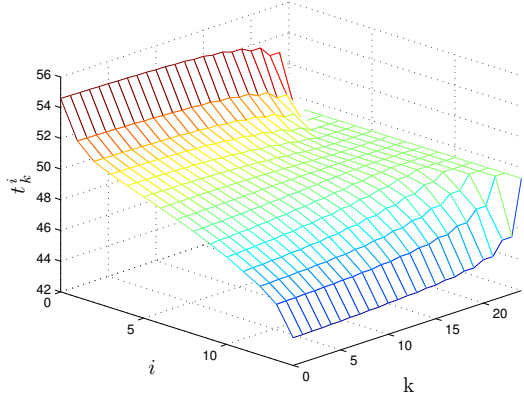}
\vspace{-5pt}
\caption{Thresholds as a function of time and state for a discretized truncated log-normal distribution (and i.i.d prices) with mean $\overline{\lambda}_{k}=49$ and $\sigma_{k}=9$ for all $k$, $\overline{v}=1$, n=$15$, $N=24$, $\hat{t}=\overline{\lambda}_{k}$ and no storage penalties. 
}
\label{threshpic}
\end{center}
\vspace{-10pt}
\end{figure}
\vspace{-3pt}

\begin{remark}
The form of the optimal policy shows that if we start with an empty storage ($s_0=0$), then the storage state $s_k$ will only take integer multiples of $\overline{v}$ since $v^\ast_{k}\in \{-\overline{v},0,\overline{v}\}$ for all $k$. When $s_0\neq0$, the storage state will fall on the grid of integer multiples of $\overline{v}$ immediately after the first time that $v^\ast_{k}=i\overline{v}-s_{k}$ or $v^\ast_{k}=(i+1)\overline{v}-s_{k}$, and we will have $v^\ast_{k}\in \{-\overline{v},0,\overline{v}\}$ for the remainder of the time horizon. This conclusion holds for the infinite-horizon case as well and can help simplifying the policy and analysis by focusing on a finite state system.
\end{remark}
\begin{remark}
The upperbound $\overline s$ on the storage capacity is enforced by choosing $h^{i}_{k}$ in (\ref{penal}) sufficiently large (i.e. $h^{i}_{k}>\lambda^{\max}_{k}$) for $i \geq n$, so that it would never be optimal to store energy beyond $\overline s$. It can be verified that the thresholds and consequently the optimal policy are invariant with respect to the choice of $h^{i}_{k}$ for $i\geq n$ as long as $h^{i}_{k}>\lambda^{\max}_{k}$.
\end{remark}
\smallbreak
\noindent{\bf Definition 2.} We define the economic value of storage, or simply the value of storage, as the negative of the cost of the optimal value of problem (\ref{storage}) divided by the number of stages ($N$), and we denote it by $\mathcal{V}$ for the finite-horizon case and by $\mathcal{V}_{\infty}$ for the infinite-horizon case. Therefore, $\mathcal{V}= -V_{0}\left( s_{0}\right)/N$. For instance, if $s_0=0$ (i.e. the consumer starts with an empty storage), it then follows from (\ref{valuefunc}) that for the finite-horizon case:
%\vspace{-3pt}
\begin{equation} \label{x0earnings}
\mathcal{V}=-V_{0}(0)=-e^{0}_{0}/N.
%\vspace{-4pt}
\end{equation}
Thus, $\mathcal{V}$ can be computed using the recursive equations in (\ref{e_k}) for the finite-horizon case. 
\subsection{The finite-horizon value of storage under empirical price distribution from real-time wholesale electricity markets} 
Herein, to examine our results and verify the validity of our assumptions, we test the optimal policy against actual market price data. We apply the finite-horizon optimal policy (\ref{optpol}) to real-time wholesale market data taken from the \cite{PJM} and the \cite{isone} to examine the competitive ratio (CR) of our policy when applied to real-world price data. This would allow us to assess whether the assumption of independent prices is actually reasonable for the storage problem. A relatively high CR would not necessarily suggest that prices are independent; rather, it could mean that the sensitivity of the optimal policy to the correlations, if any, in real-world prices is low. If the CR obtained form the assumption of independent prices is relatively high, we can suggest that the additional information that the current price provides on future prices (as a conditional distribution) has little value. 

The setup of the experiment is as follows. We take the actual data for hourly energy prices (for $16$ hours of each day, from $8$ a.m. to $12$ a.m.) for two different months (December 2010 and July 2011) from the \cite{PJM}, and May and November 2011 from \cite{isone}. The choice of these dates and times was arbitrary. For the purpose of these simulations, we set $\hat{t}$, the salvage value, equal to $\overline{\lambda}$, the empirical mean. To perform the simulations, we first find the empirical distribution of the data for each case. In these computations, we use the same price distribution for all $k$; i.e., we take all the price data for the entire month, and use this data to find the empirical distribution of prices in that month and compute the thresholds for the optimal policy for all hours from the same empirical distribution. In other words, we assume that prices are independently and identically distributed (i.i.d). Then, at the beginning of each stage we reveal the actual price of that stage to the optimal policy and record the decision. This gives us the profit resulting from the optimal policy. Next, for the purpose of comparison, we assume the entire price sequence is perfectly known a priori and compute the value that results from deterministically and omnisciently maximizing the profit against the materialized prices. This deterministic value is the absolute best that an omniscient agent could have done, and we will refer to its corresponding deterministic policy as the omniscient policy. We then find the CR by computing the ratio of the value obtained from the optimal policy to the value obtained from the omniscient policy. The CR gives us a measure of how well our optimal policy has performed. 

For each month, there is a number of days in which the average price is well above the average price of the entire month. These days are outliers in the sense that the empirical distribution is way off for modeling their price sequence. We go about this issue by performing three sets of simulations. For the first set of simulations, we only select those days in which the average price is within one standard deviation of the average price of all the selected days in the month. In the second set of simulations, we select those days in which the average price is within about 1.5 times the standard deviation of the average price of all the selected days in the month. Then, in the third set of simulations, we take all the days of the month into account, even the outliers. We record the CR for each day, and report the average CR for each month and each case in Table \ref{tbl1}. Comparing the results in Table \ref{tbl1} reveals how much these outlier days affect the CR. Note that a ramp constraint of $\overline{v}=1$ and a storage capacity of $\overline{s}=10$ is used in all the simulations, and it is assumed that there are no penalties on storage.

\begin{table}[tp]%
\caption{Competitive Ratio for each set of simulations}
\label{tbl1}\centering%
\begin{tabular}{lccc}
\hline\hline
Source/Month  & Competitive Ratio & \# of qualifying days & Range from the mean	\\ \hline 
PJM / July	 & 0.86&  27 &$1\times \sigma$		\\
 & 0.84& 	28&1.5$\times \sigma$	\\
 & 0.81&	All	&All	\\ \hline 
PJM / December &0.88	&  22	&$1\times \sigma$\\
 &0.87	& 25&$1.5\times \sigma$\\
&0.77	&	All&All \\ \hline 
ISO NE/ May  	& 0.72  & 25&1$\times \sigma	$	\\
 	& 0.71	& All&$1.5\times \sigma$\\
	& 0.71 &	All &All	\\ \hline 
ISO NE/ November	& 0.89	& 29&$1\times \sigma$	\\
	& 0.88&  All	&$1.5\times \sigma$\\
	& 0.88	&All	&All	\\ \hline \hline  
\end{tabular}
\vspace{-5pt}
\end{table}

\renewcommand{\thetable}{\arabic{table}}

\begin{figure}
[ptbh]
\begin{center}
\vspace{-6pt}
\includegraphics[scale=0.41]{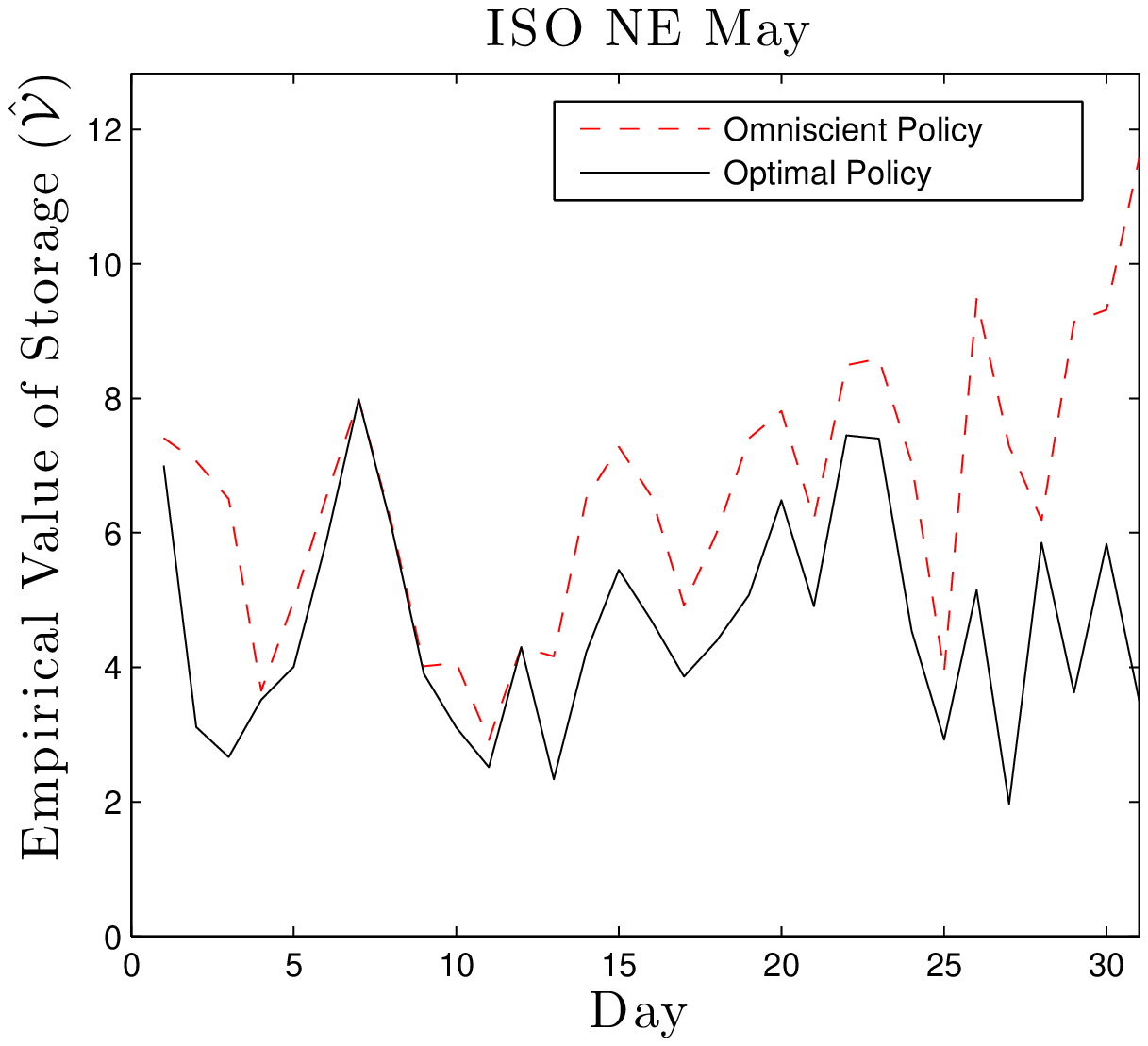}
\includegraphics[scale=0.40]{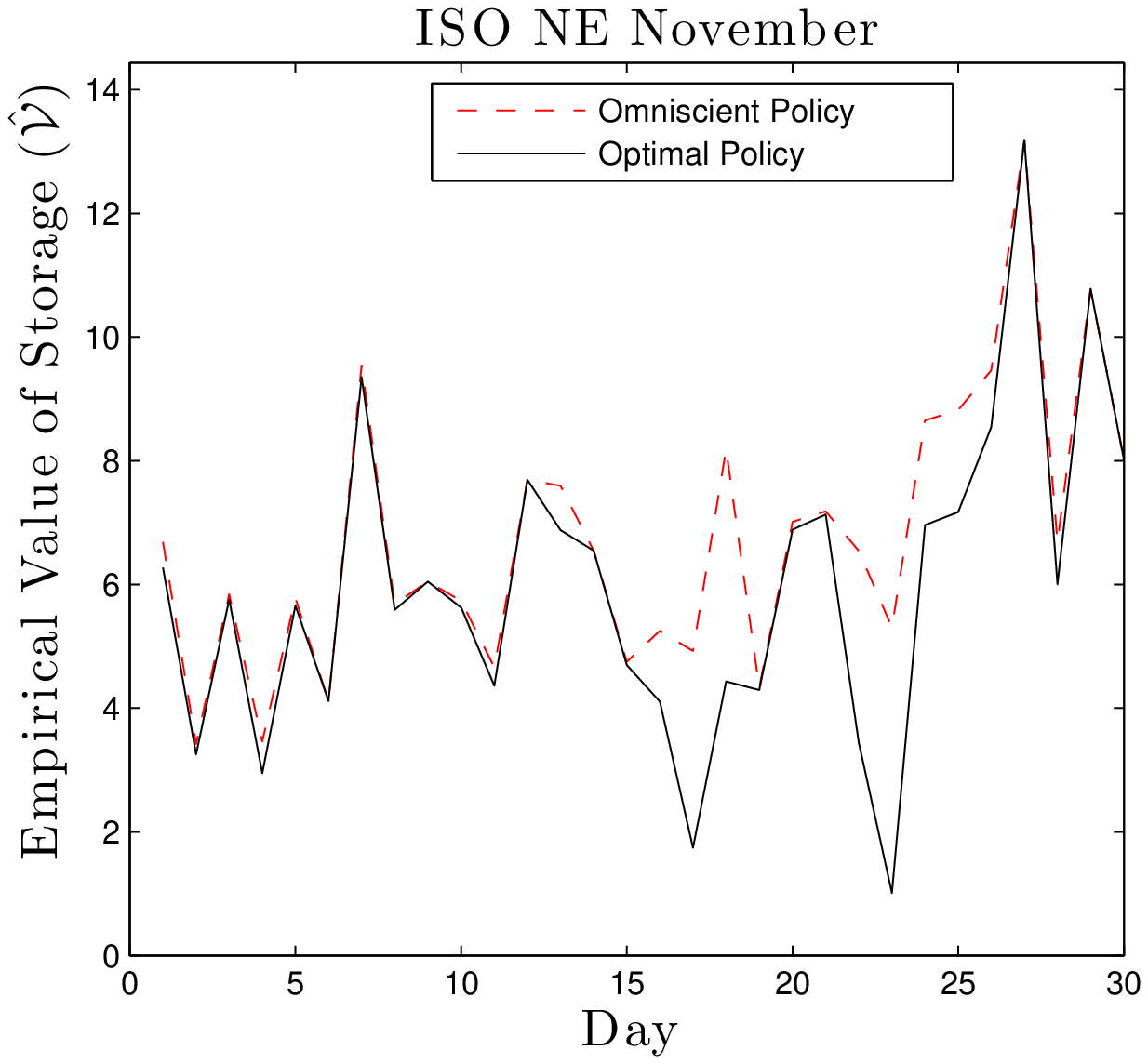}
\includegraphics[scale=0.40]{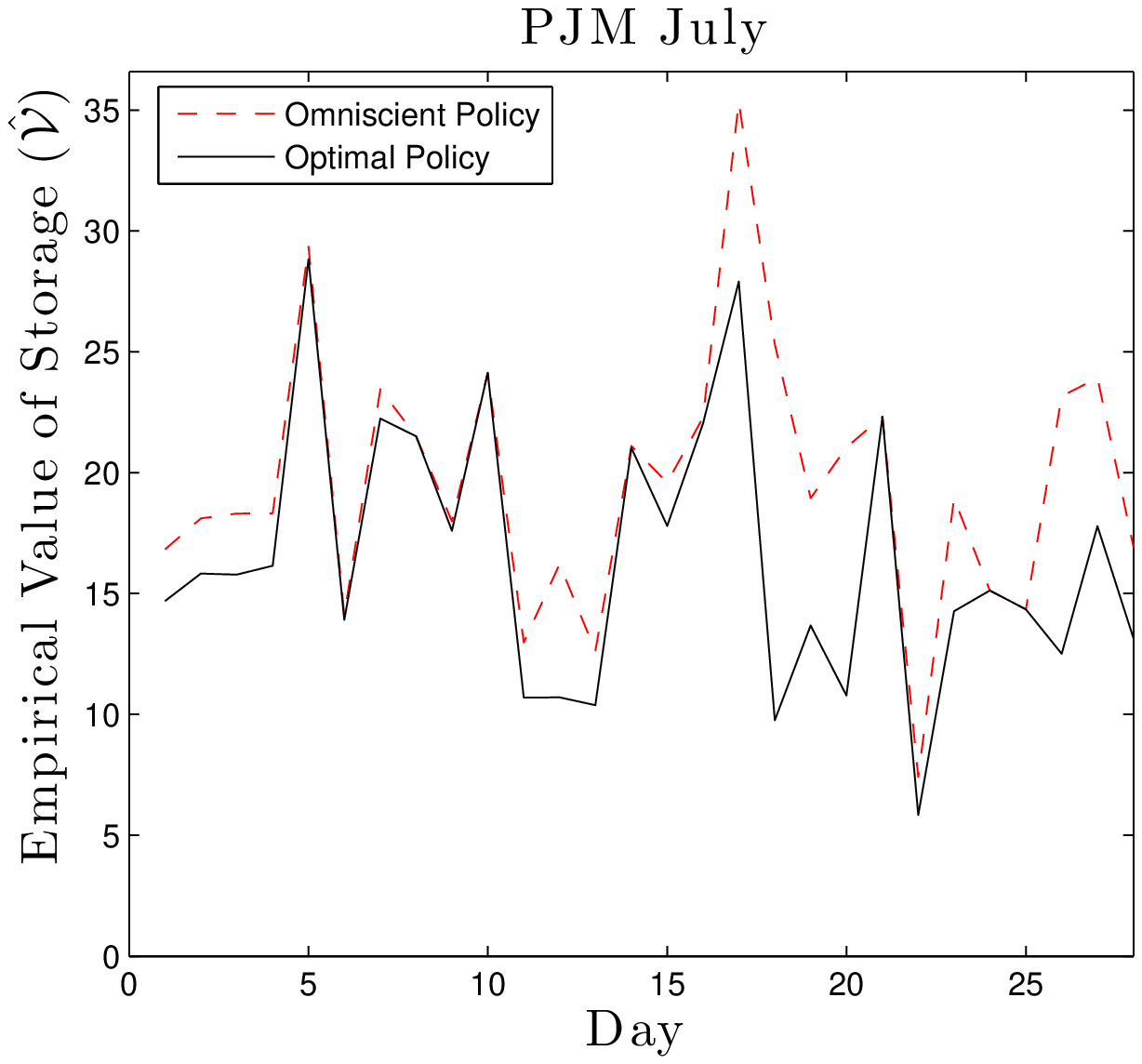}
\includegraphics[scale=0.41]{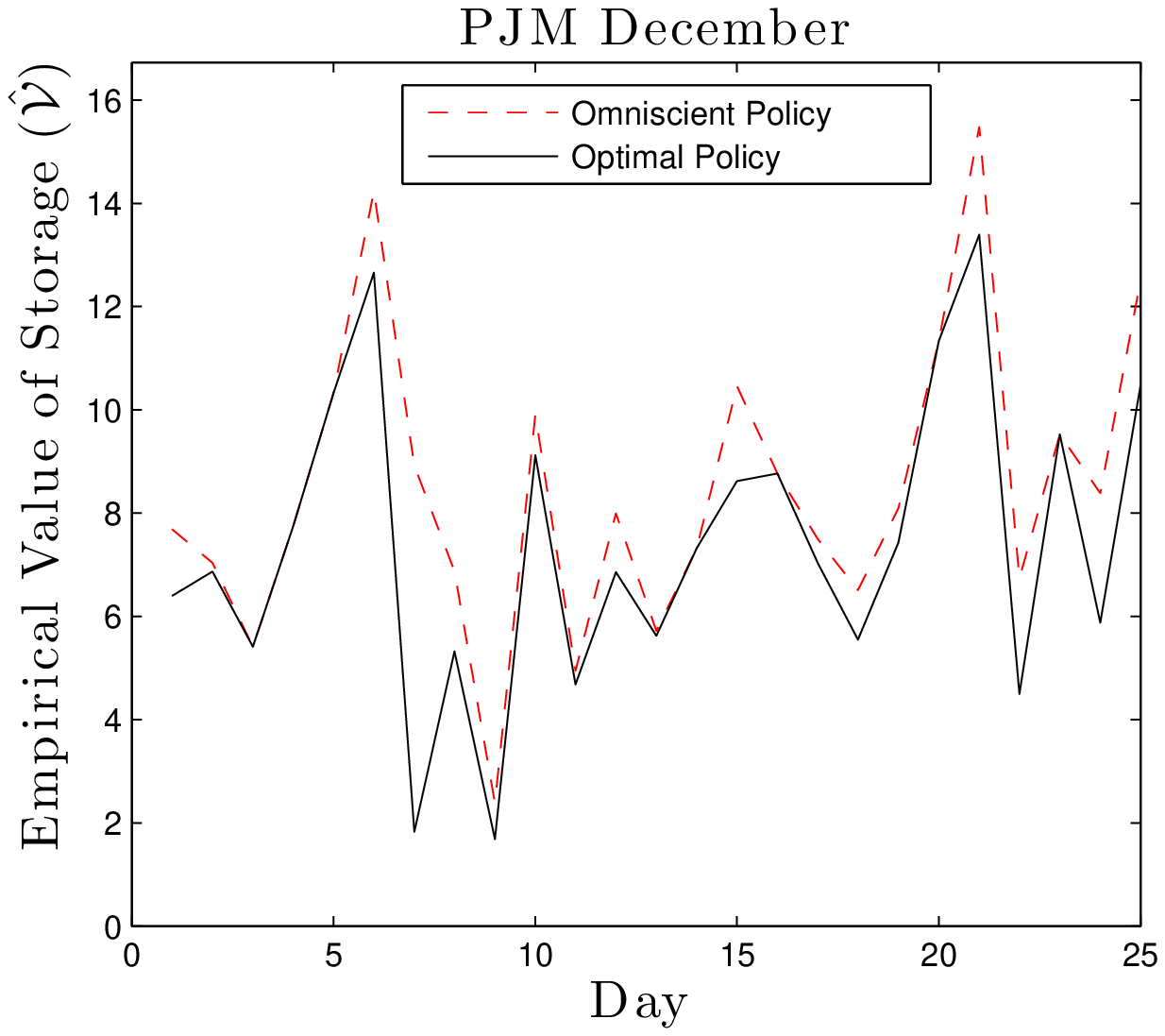}
%\vspace{-10pt}
\caption{The empirical value of storage for each day of the month in the second set of simulations, obtained from applying the optimal policy (solid line) and the omniscient policy (dashed line)}
\label{valcomp}
\end{center}
\vspace{-5pt}
\end{figure}

Figure \ref{valcomp} shows the plots of the empirical value of storage for each day of the month in the second set of simulations (i.e for those days whose average price is within 1.5 times the standard deviation of the average price of all the selected days in the month). 

As it can be seen in Figure \ref{valcomp}, we have nearly perfect matching in some days, while in some other days there is discrepancy. This discrepancy is mainly due to two reasons. First, it is because we have multiple spikes in some days with only a few (or no) prices that are below our buying thresholds. So, even though there are ample opportunities for arbitrage, even the lowest of these spiky prices is not within the normal range. In Table \ref{tbl1}, this effect can be observed for the month of December in PJM, in which the removal of those days with high average prices improved the CR from $0.77$ to $0.88$. The second reason is that in some days all the prices are almost in the same range compared to our thresholds (i.e. they are either mostly above our buying thresholds or mostly below our selling thresholds). In other words, the low CR is just an outcome of an undesirable sequence of prices (sample path). So, even though for the deterministic case with the price sequence known a priori it is possible to take advantage of these small price differentials, our thresholds are unable to capture these opportunities. This effect can be observed in the month of May in Table \ref{tbl1}, in which the removal of those days with high average prices did not really improve the CR; this is because the average price in the days with a relatively flat price sequence is not necessarily considerably higher than the month's average.

So far, we have been computing the thresholds for the optimal policy using the empirical price distributions from the prices in that same month. This can be taken as a proxy for the sensitivity of the value of storage to correlations in the actual prices, given that our optimal policy assumes independent prices. Although there is no benchmark to compare with, the CR obtained from applying the optimal policy seems reasonably high, and we suggest that the sensitivity of the value of storage to the assumption of independent prices is low. Also, although we do not do as well on those days with higher than normal, or very flat price profiles, it does not appear that a Markovian or Martingale assumption on prices %
could do any better, because these prices are really outliers and do not seem to follow a structured stochastic pattern that can be learned from past data. However, this needs to be substantiated by further studies and systemic experiments.

In the next set of experiments, instead of computing the empirical distribution from the data of that same month, we compute the empirical distributions from the price data of the past 30 days, the past 20 days, and the past 10 days, respectively. The results are shown in Table \ref{tbl4}.

\begin{table}[h]  \caption{Competitive ratios using the empirical distribution from the past 30, 20, and 10 days} %title of the table
\centering \begin{tabular}{c c c c c} \hline\hline Source / Month & Past 30 days  & Past 20 days & Past 10 days &\# of selected days \\ [0.5ex] \hline 
 PJM / July &  0.68& 0.70 &0.75& All \\	% Entering row contents 
PJM / December  &0.55&0.65&0.67& All\\ 
ISO NE/ May  	& 0.69 &0.70&0.68& All\\
ISO NE/ November 	& 0.91&0.89&0.89& All \\[1ex] 
\hline \end{tabular} \label{tbl4} \end{table}

An interesting observation is that for both months of May and November in ISO NE, using the empirical distribution from the past 20 days gives almost the same CR as using the empirical distribution from May and November themselves (as reported in Table \ref{tbl1}). However, comparing the competitive ratios shown in Table \ref{tbl4} for both months in PJM with the results shown in Table \ref{tbl1}, we observe that for PJM, using the empirical distribution from historical data does not do as well as using the empirical distribution from that month itself. %T

%\vspace{-3pt}
%%%%%%%%%%%%%%%%%%%%%%%%%%%%%%%%%%%%%%%%%%%%%%%%%%%%%%%%%%%%%%%%%%%%%%%%%%%%%%%%
%%%%%%%%%%%%%%%%%%%%%%%%%%%%%%%%%%%%%%%%%%%%%%%%%%%%%%%%%%%%%%%%%%%%%%%%%%%%%%%%
%%%%%%%%%%%%%%%%%%%%%%%%%%%%%%%%%%%%%%%%%%%%%%%%%%%%%%%%%%%%%%%%%%%%%%%%%%%%%%%%
\section{Characterization of the Economic Value of Storage} \label{valchar}
\subsection{Analytical upperbound on the value of storage} \label{analyticbound}
In this subsection, we will derive a bound on the long-term economic value of ramp-constrained storage. Herein, for the purpose of obtaining the bound, we further assume that $h^{i}_{k}=0$, for all $i\leq n-1$ and $k\leq N$. This assumption is consistent with our objective of finding an upperbound on the value of storage.
\smallbreak
\noindent{\bf Definition 3.} Given a control policy $\pi_k:[0,\infty)^2\mapsto [-\bar{v},\bar{v}]$, and starting from an arbitrary initial state $s$, the infinite-horizon average cost per stage associated with problem \eqref{storage} is defined as
\begin{equation}
\begin{aligned}
%\min~~
& \gamma_\pi \overset{\text{def}}=\underset{N\rightarrow \infty}{\lim}
& & \!\!\frac{1}{N}~  \mathbf{E}\left[{\displaystyle\sum\limits_{k=0}^{N-1}} \lambda_{k}v_k ~|~ s_0=s\right],\\
\end{aligned}
%J_\pi(s)=
\label{eq:avgcost}
\end{equation}
where $v_k=\pi_k(x_k,\lambda_k).$ We will refer to the problem of optimization of $\gamma_\pi$ over all feasible stationary policies as the infinite-horizon storage management problem. We will denote the associated optimal cost by $\gamma^\ast$ and refer to 
$\mathcal{V}_{\infty}\overset{\text{def}}=-\gamma^\ast$ as the long-term expected economic value of storage.
%\end{definition*}
\begin{remark}
It is standard to show that the optimal average cost is independent of the initial state $s_0$. Moreover, if the relative value iteration for the infinite-horizon storage problem converges to some differential cost function $H^{\ast}(s)$, then it is necessary for $H^{\ast}(s)$ and the optimal average cost per stage $\gamma^{\ast}$ to satisfy the Bellman equation (see, for instance, \citealt{Berts}):
\begin{equation}
H^{\ast}(s)=\mathbf{E} \left[\min_{v \in [\max(-s,-\overline{v}),\min(\overline{v},\overline{s}-s)]} \lambda v +H^{\ast}(s+v)\right]-\gamma^{\ast}.  \label{valiter}
\end{equation}
\label{rem:initC}
\vspace{-10pt}
\end{remark}
\begin{theorem}
Consider the infinite-horizon storage management problem. Suppose that the support of the price distribution function at all stages lies within an interval $[\lambda_{\min},\lambda_{\max}] \subseteq [0,\infty)$. All else held constant, the maximum over all possible distributions, of the long-term economic value of storage is given by 
\begin{equation}
\mathcal{V}_{\infty}=-\gamma^{\ast}=\frac{\overline{v}(\lambda_{\max}-\lambda_{\min})}{2}\frac{n}{n+1}=\frac{(\lambda_{\max}-\lambda_{\min})}{2}\frac{\overline{s}~\overline{v}}{\overline{s}+\overline{v}},
\label{2ptgamma}
\end{equation}
and is attained when the prices are sampled from a two-point symmetric distribution with nonzero probability masses placed at the endpoints of the fixed support:\vspace{0.1in} 
\begin{equation}
P_{\Lambda}(\lambda)=\left\{ 
\begin{array}
[c]{lcl}
1/2 &\text{if}& \lambda = \lambda_{\min} \\[0.05in]
1/2 &\text{if}& \lambda = \lambda_{\max} \\[0.05in]
0 & \text{otherwise}
\end{array}\label{2ptdist}
\right. 
\end{equation}
Furthermore, the corresponding differential cost function satisfying the Bellman equation is given by the following piecewise linear convex function: 
\begin{equation}
 H^{\ast}(s)=-\frac{(i+1)\lambda_{\min}+(n-i)\lambda_{\max}}{n+1}s-\frac{i(i+1)(\lambda_{\max}-\lambda_{\min})\overline{v}}{2(n+1)}, \label{diffcost}
\end{equation} 
where $s\in\left[i\overline
{v},(i+1)\overline{v}\right)$ and $i\in \{0,\cdots,n-1\}$, and for the special case of $s=\overline{s},$ we take $i=n-1$.

\noindent Finally, for any two-point distribution with PMF
\begin{equation}
P_{\Lambda}(\lambda)=\left\{ 
\begin{array}
[c]{lcl}
a &\text{if}& \lambda = \lambda_{\max} \\[0.05in]
1-a &\text{if}& \lambda = \lambda_{\min} \\[0.05in]
0 & \text{otherwise}
\end{array}\label{gen2ptdist}
\right. 
\end{equation}
we have
\begin{equation}
\mathcal{V}_{\infty}=-\gamma^{\ast}=\overline{v}(\lambda_{\max}-\lambda_{min})\frac{b(1+b+\cdots+b^{n-1})}{(b+1)(1+b+\cdots+b^{n})} \label{gengamma}
\end{equation}
where $b=(1-a)/a$.
\label{2ptsymm}
\end{theorem}
\proof{Proof.}
Please see the Appendix.
\Halmos
\endproof

\smallbreak
If in addition to the support of the price distribution we also fix the mean of the distribution, we can obtain a tighter bound as stated in Corollary \ref{meancoroll}:

\begin{corollary}
Suppose that in addition to fixing the support of the price distribution function, we also fix the mean of the price distribution to $\mu \in (\lambda_{\min},\lambda_{\max})$. All else held constant, the maximum over all possible distributions, of the long-term expected value of storage is attained when the prices are sampled from a two-point distribution with the following PMF:\vspace{0.1in} 
%P_\ 
\begin{equation*}
P_{\Lambda}(\lambda)=\left\{ 
\begin{array}
[c]{lcl}
\frac{\mu-\lambda_{\min}}{\lambda_{\max}-\lambda_{\min}} &\text{if}& \lambda = \lambda_{\max} \vspace{5pt}\\
\frac{\lambda_{\max}-\mu}{\lambda_{\max}-\lambda_{\min}}  &\text{if}& \lambda = \lambda_{\min}\vspace{5pt} \\
0 & \text{otherwise}
\end{array}%\
\right. 
\end{equation*}
The corresponding long-term value of storage is obtained by plugging $b=(\lambda_{\max}-\mu)/(\mu-\lambda_{\min})$ into \eqref{gengamma}.
\label{meancoroll}
\end{corollary}

\proof{Proof.}
Please see the Appendix.
\Halmos
\endproof

\begin{remark}
The $n/(n+1)$ scaling in the optimal average cost per stage implies that $90\%$ of the maximum possible value of storage is achieved when the storage capacity is only 9 times the ramp constraint. Note that this was obtained for an extreme distribution. As we will see in the remainder of this section, for less extreme distributions with smaller variance over the support, the value of storage saturates even more quickly. This includes empirical distributions obtained from electricity market data. Furthermore, aging, dissipation, and non-ideal charging and discharging factors further reduce the value of storage. 
\end{remark}

\subsection{Computational Experiments for Characterization of the Economic Value of Storage}
%%%%%%%%%%%%%%%%%%%%%%%%%%%%%%%%%%%%%%%%%%%%%%%%%%%%%%%%%%%%%%%%%%%%%%%%%%%%%%%%
%%%%%%%%%%%%%%%%%%%%%%%%%%%%%%%%%%%%%%%%%%%%%%%%%%%%%%%%%%%%%%%%%%%%%%%%%%%%%%%%
%%%%%%%%%%%%%%%%%%%%%%%%%%%%%%%%%%%%%%%%%%%%%%%%%%%%%%%%%%%%%%%%%%%%%%%%%%%%%%%%
In this subsection, we employ numerical computations to characterize the economic value of the proposed model of storage over a finite time-horizon and highlight the effects of ramp constraints and price volatility on the value of storage.
\label{sec:earnings}

\begin{figure*}
[ptb]
\begin{center}
\vspace{-4pt}
\includegraphics[scale=0.37]{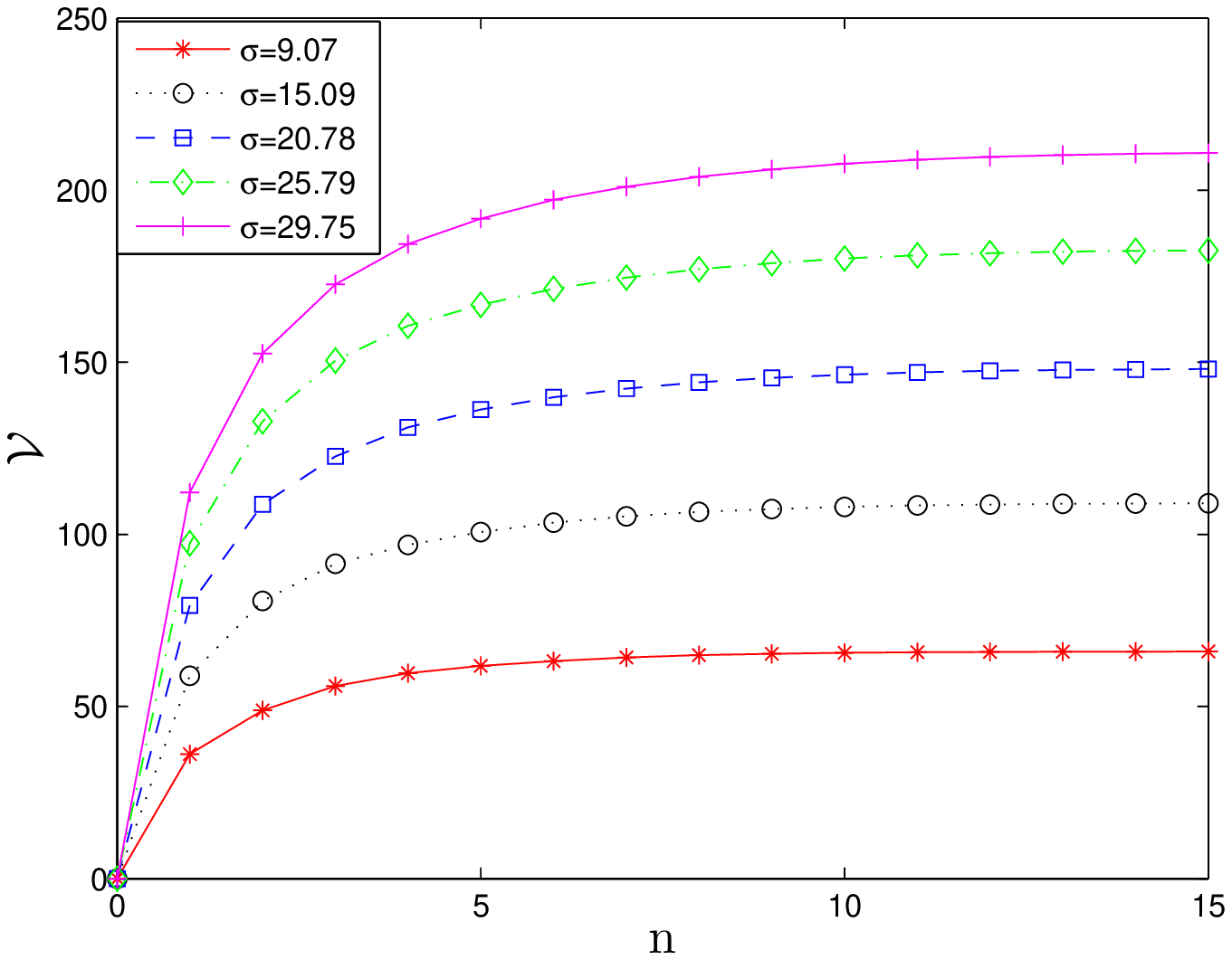}%{n.eps}
\includegraphics[scale=0.37]{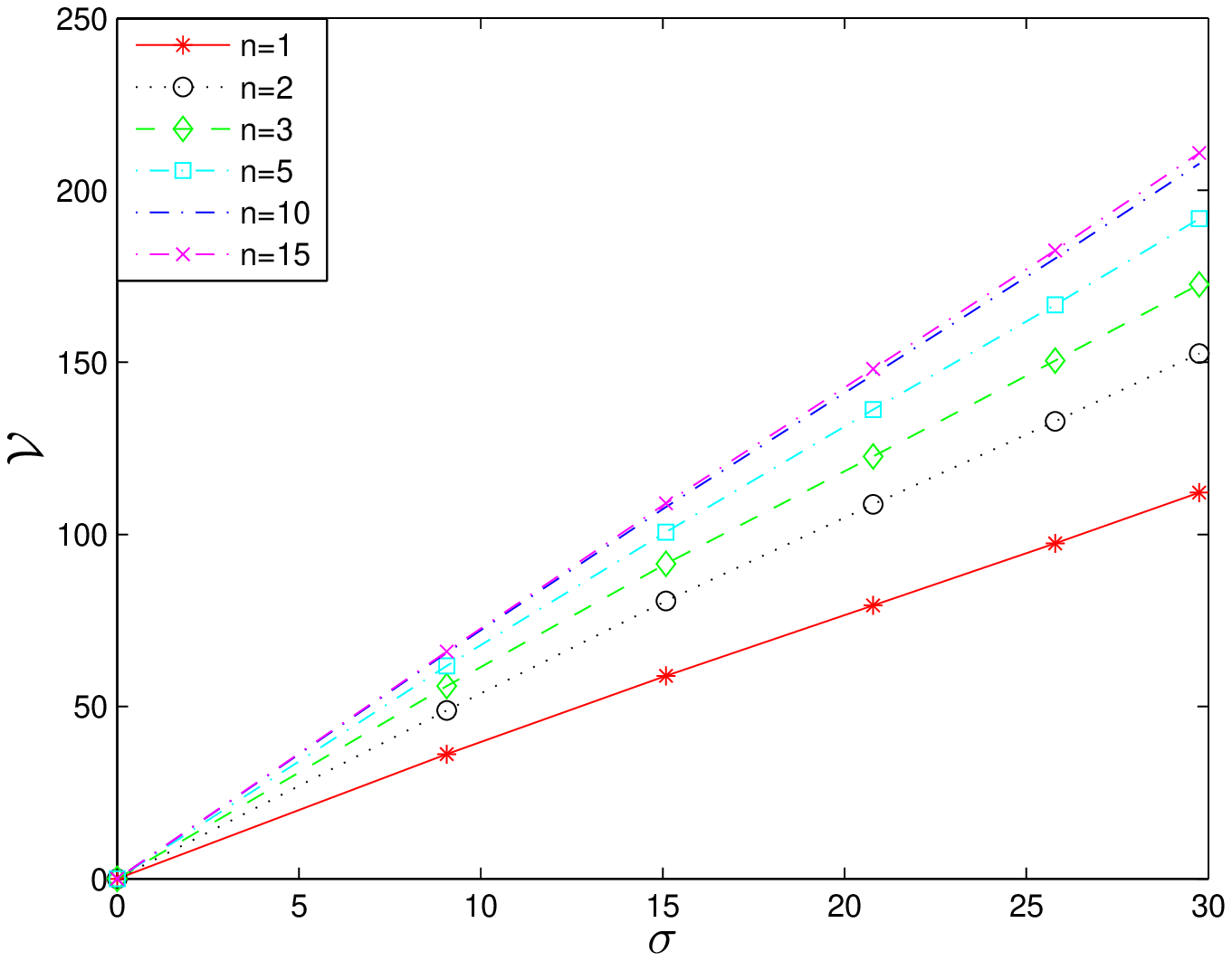}%{sig.eps}
\vspace{-7pt}
\caption{\hspace{-3pt}$\mathcal{V}$ vs. $n$ (left) and $\sigma$ (right) for a few samples, using a discretized truncated log-normal distribution.}
\label{fig:scenario-a}
\end{center}
\vspace{-7pt}
\end{figure*}

\begin{figure*}
[ptb]
\begin{center}
\vspace{-6pt}
\includegraphics[scale=0.37]{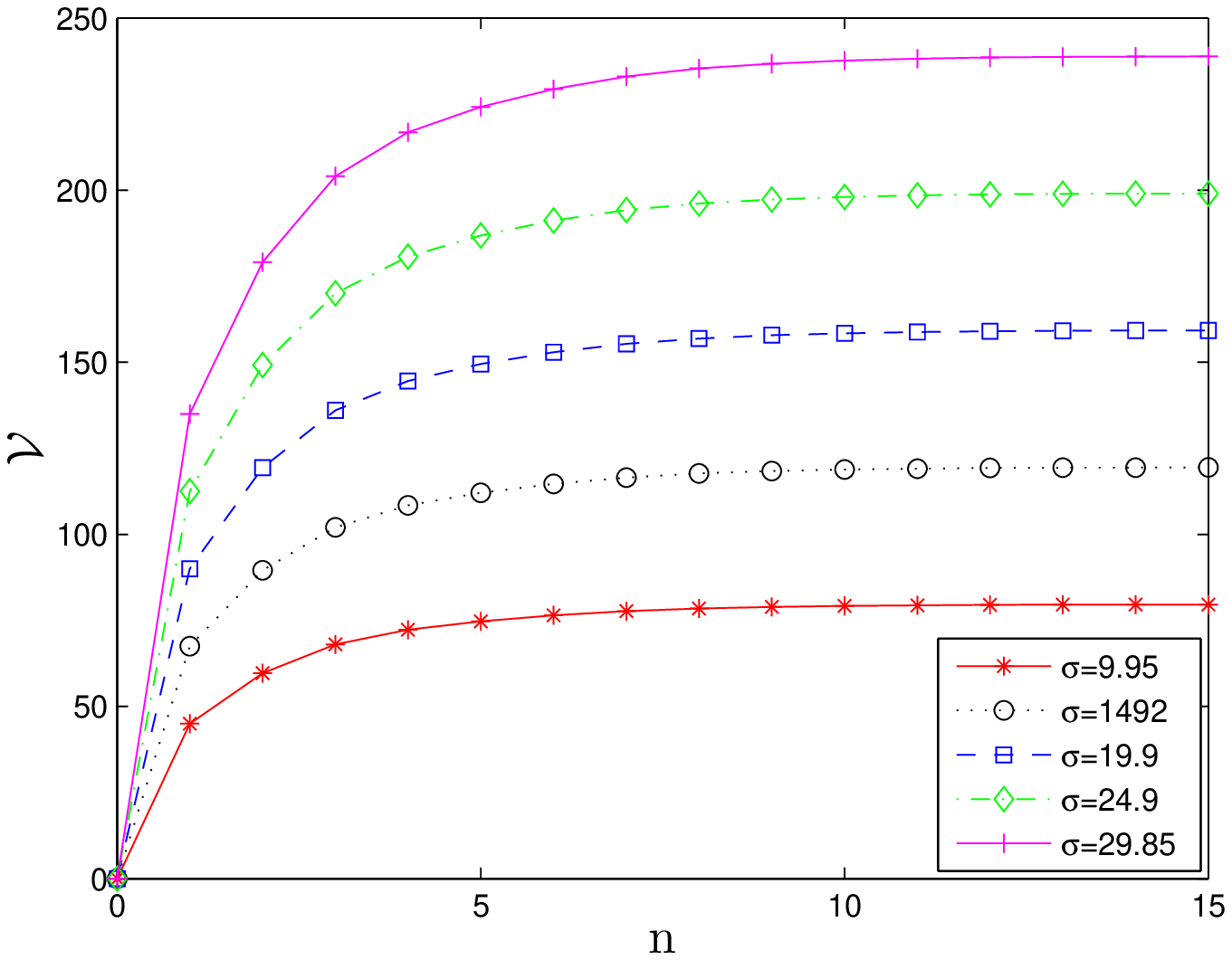}%{HiLo_n.eps}
\includegraphics[scale=0.37]{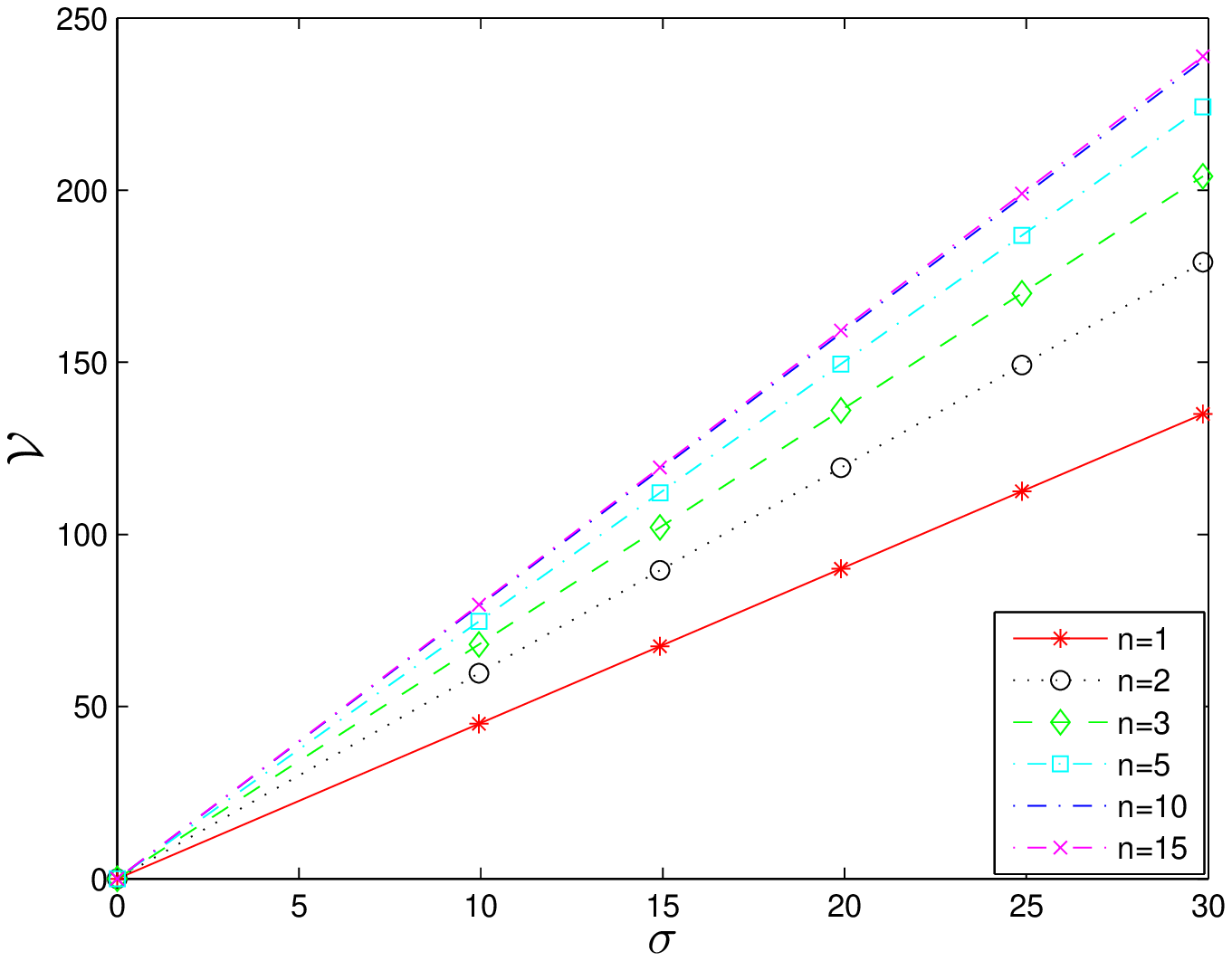}%{HiLo_sig.eps}
%\vspace{-10pt}
\caption{$\mathcal{V}$ vs. $n$ (left) and $\sigma$ (right) for a few samples, using a discretized uniform distribution.}
\label{fig:scenario-aa}
\end{center}
\vspace{-8pt}
\end{figure*}

\smallbreak
\noindent Herein, we consider the following classes of distributions:
\begin{list}{\labelitemi}{\leftmargin=1em}
\item Discretized truncated log-normal distribution, with fixed mean $\overline{\lambda}=50$, 
\item Discretized uniform distribution, with fixed mean $\overline{\lambda}=50$.
\end{list}
The reason for choosing the log-normal distribution is that the empirical distributions from ISO New England and PJM qualitatively resemble a log-normal distribution, at least for the cases that we tested. The choice of the mean ($\overline{\lambda}=50$) is also a realistic choice for average hourly energy prices in markets such as PJM and ISO NE. For the purpose of these computations, we use the same price distribution for all $k$ (i.e. we assume that prices are independently and identically distributed). Throughout this section we assume that $s_0=0$, which means that the consumer starts with an empty storage.
For each of these distributions, we fix all quantities in our model other than $\sigma$, the standard deviation of price distribution, and $n$, the ratio of storage capacity $\overline{s}$ to physical ramp constraint of storage $\overline{v}$. We vary $n= \overline{s}/\overline{v}$ by fixing $\overline{v}$ and changing $\overline{s}$. Using the fixed quantities $N=24$, $\overline{v}=10$, and $\overline{\lambda}=50$, we examine how $\mathcal{V}$ varies as a function of $\sigma$ and $n$ for the case of no storage penalties. For the purpose of these simulations, we set $\hat{t}$ equal to the mean of the price distribution. %
%\vspace{-3pt}
Herein, we set $h^{i}_{k}=0$, for all $i\leq n-1$ and $k\leq N$, so that there is no penalty on storing energy up to capacity. Then, for a fixed time horizon, we examine how $\mathcal{V}$ varies with $\sigma$ and $n$, for each of the following price distributions:\\
\emph{Discretized truncated log-normal distribution:} Figure \ref{fig:scenario-a} illustrates how $\mathcal{V}$ changes with $\sigma$ and $n$, for the discretized truncated log-normal distribution. The plots show that the value of storage increases linearly with $\sigma$. As one would expect, the value of storage also increases as the storage capacity increases. However, it is interesting to note that for a fixed standard deviation, the value of storage saturates fairly quickly as a function of $n$. Hence, for a given time horizon, a fixed ramp constraint, and a fixed $\sigma$, there exists a certain range for capacity beyond which the value of storage will no longer change noticeably. Also, the optimal storage capacity increases with price volatility. \\
\emph{Discretized uniform distribution:} 
As can be seen in Figure \ref{fig:scenario-aa}, saturation of the value of storage occurs at about the same value as in the log-normal case. However, for certain extreme distributions, such as an asymmetric 2-point distribution, saturation can occur more quickly. Note also that the value of storage is a linear function of the standard deviation, just like the log-normal case.

%\vspace{-6pt}
%\vspace{.1in}

One interesting observation in these results is that in the presence of ramp constraints, several distributed storage systems would be more profitable than one large storage system of equal ramp constraint and aggregate capacity, due to the quick saturation of $\mathcal{V}$ as $n$ increases. Although this observation is based on the assumption that the ramp constraint and capacity are independent, this assumption might actually be valid for the case of distribution grids in which the ramp constraint is imposed by the power lines. Another interesting observation is that although the shape of the plots look quite similar for both distributions, the value of the uniform distribution is higher than that of the log-normal. This makes intuitive sense because with the uniform distribution, we have as many opportunities for buying at a low price as we have for selling at a high price, while for the log-normal case, most of the probability mass is centered around the mode, which creates fewer opportunities for arbitrage.

\paragraph{\bf Using wholesale market data from ISO New England:} For the purpose of comparison, we repeat the experiment using the data for hourly energy prices from three different months (April, June, and October of 2011) obtained from \cite{isone}. For each month, the hourly data (from all 24 hours of all days of the month) have been used to find the empirical distribution of prices for that particular month. For the purpose of these computations, we compute the empirical distribution over all stages, and we use this empirical distribution for all the stages (i.e. we assume that prices are i.i.d). Figure \ref{isonesat} shows how $\mathcal{V}$ varies with $n$ in each month without storage penalties:
\vspace{-3pt}
\begin{figure}
[ptbh]
\begin{center}
\vspace{-15pt}
\includegraphics[scale=0.375]{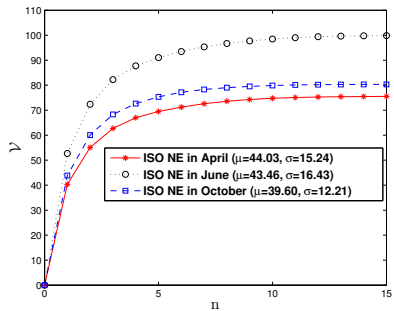}
%\
\vspace{-9pt}
\caption{$\mathcal{V}$ vs. $n$ for 3 different months, without storage penalties}
\label{isonesat}
\end{center}
\vspace{-15pt}
\end{figure}

Note in Figure \ref{isonesat} that similar to the results of the distributions used in Section \ref{sec:earnings}, the value of storage nearly saturates when $n$ is about 6. An interesting observation is that although both the mean and the standard deviation of the empirical price distribution in October are lower than those of the empirical price distribution in April, the value of storage is higher in October. This observation reiterates the importance of the price distribution for the value of storage, even when its mean and standard deviation are somewhat lower, though many other factors such as non-stationarity and correlation in the data could skew the results.

%\vspace{-10pt}
%%%%%%%%%%%%%%%%%%%%%%%%%%%%%%%%%%%%%%%%%%%%%%%%%%%%%%%%%%%%%%%%%%%%%%%%%%%%%%%%
\section{Implications of the Optimal Policy}\label{sec:policyimp}
Herein, we will focus on understanding the price-dependence and state-dependence of the optimal policy, and the impact of the state-price interplay on the storage response. We first average out state-dependence of storage response and focus our attention on ``average" price-responsiveness only, which allows us to address the expected price elasticity of demand induced by storage in the first subsection. Then, in the next two subsections, we focus on the interplay between state-dependence and price-dependence of storage response, and study potential impacts of the state-and-price-dependent response from market-based storage on the required reserves.
\subsection{Average Price Elasticity of Demand Induced by Storage}
%%%%%%%%%%%%%%%%%%%%%%%%%%%%%%%%%%%%%%%%%%%%%%%%%%%%%%%%%%%%%%%%%%%%%%%%%%%%%%%%
%%%%%%%%%%%%%%%%%%%%%%%%%%%%%%%%%%%%%%%%%%%%%%%%%%%%%%%%%%%%%%%%%%%%%%%%%%%%%%%%

\label{sec:elasticity}
%\vspace{-2pt}
In this subsection, in order to study the implications of the optimal policy on price elasticity of demand (PED) in an electricity market, we present a computational framework for studying the average PED in a simulated electricity market. In our dynamic model, the storage response depends on price, stage, state, time-horizon, storage capacity, and ramp constraint. The term ``average PED" reflects the fact that the dependence of storage response on the stage and the internal state of storage have been averaged out. We assume that there is a fixed time horizon $N$, and eliminate state-dependence by taking expectations. In particular, we define:
\vspace{-5pt}
\begin{equation*}
v(k,\lambda) =\mathbf{E}_{s_0,s_k}\left[  v^{\ast}_{k}|\lambda_{k}=\lambda\right].
\vspace{-5pt}
\end{equation*}
In order to eliminate stage-dependence, we think of the storage response-measuring observer as sampling a random time $\tau$ uniformly over $\{0,\cdots,N\}$. By averaging over this randomness, we maintain dependence on price alone:
\vspace{-5pt}
\begin{equation*}
v_\text{avg}(\lambda) =  \mathbf{E}_{\tau} \left[v(\tau,\lambda) \right],
\vspace{-5pt}
\end{equation*}
which is captured in the simulations by clustering real-time prices, and averaging over each cluster. 
%\vspace{-5pt}
%%%%%%%%%%%%%%%%%%%%%%%%%%%%%%%%%%%%%%%%%%%%%%%%%%%%%%%%%%%%%%%%%%%%%%%%%%%%%%%%
%%%%%%%%%%%%%%%%%%%%%%%%%%%%%%%%%%%%%%%%%%%%%%%%%%%%%%%%%%%%%%%%%%%%%%%%%%%%%%%%
%\subsection{Simulations}
%%%%%%%%%%%%%%%%%%%%%%%%%%%%%%%%%%%%%%%%%%%%%%%%%%%%%%%%%%%%%%%%%%%%%%%%%%%%%%%%
%%%%%%%%%%%%%%%%%%%%%%%%%%%%%%%%%%%%%%%%%%%%%%%%%%%%%%%%%%%%%%%%%%%%%%%%%%%%%%%%
%\vspace{-3pt}
%%%%%%%%%%%%%%%%%%%%%%%%%%%%%%%%%%%%%%%%%%%%%%%%%%%%%%%%%%%%%%%%%%%%%%%%%%%%%%%%
%
%%%%%%%%%%%%%%%%%%%%%%%%%%%%%%%%%%%%%%%%%%%%%%%%%%%%%%%%%%%%%%%%%%%%%%%%%%%%%%%%

In these numerical simulations, we average over random instances of price and storage initial states. We set $N=288$, which corresponds to a period of 24 hours, where real-time prices are updated once in every 5 minutes. The storage system implements the optimal policy given in Theorem \ref{storageTHM}. For the purpose of these computations, we use the same price distribution for all $k$ (i.e. we assume that prices are independently and identically distributed). For generating random price sequences, we simulate a discretized truncated log-normal distribution with $\overline\lambda=52$ and $\sigma=22$, because the log-normal distribution qualitatively resembles the empirical distribution of prices from ISO NE and PJM, at least for the cases that we tested. Also, a mean of 52 and standard deviation of 22 are realistic choices for real-time energy prices in markets such as PJM and ISO NE on a day with moderate volatility. Based on our results in Section \ref{sec:earnings}, for these model parameters, a storage capacity of $\overline{s}=5\overline{v}$ is a reasonable choice for all consumers. We set $\hat{t}$ equal to $\overline\lambda$ in these simulations, just like all previous simulations. We also set $\overline{v}=10$, and $h^{i}_{k}=0$ for all $k\leq N$ and $i<n$. 

%%%%%%%%%%%%%%%%%%%%%%%%%%%%%%%%%%%%%%%%%%%%%%%%%%%%%%%%%%%%%%%%%%%%%%%%%%%%%%%%
%\su
%%%%%%%%%%%%%%%%%%%%%%%%%%%%%%%%%%%%%%%%%%%%%%%%%%%%%%%%%%%%%%%%%%%%%%%%%%%%%%%%

\begin{figure}
[ptbh]
\begin{center}
\vspace{-8pt}
\hspace{-10pt}\includegraphics[scale=0.4]{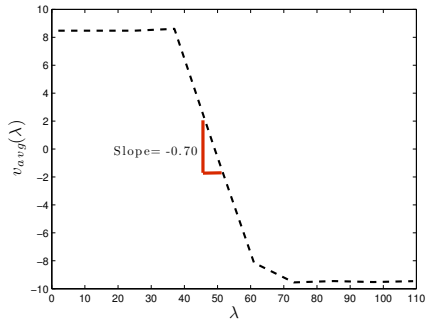}%\i
%\includegraphics[scale=0.4]{logn_aggregate_withPenalty.eps}%}
%\vspace{-11pt}
\caption{Average storage response vs. price, using the discretized truncated log-normal price distribution.}
\label{figa}
\end{center}
\vspace{-6pt}
\end{figure}

Figure \ref{figa} illustrates how the average storage response changes as a function of price.
%%%%%%%%%%%%%%%%%%%%%%%%%%%%%%%%%%%%%%%%%%%%%%%%%%%%%%%%%%%%%%%%%%%%%%%%%%%%%%%%%

As the plot suggests, the expected demand seems to be considerably more responsive to changes in the prices that fall in the mid-portion of the price range. This portion serves as a steep transition region, in which the policy quickly switches from the ``buy" policy to the ``sell" policy.

To characterize price elasticity, let us first recall the standard definition of PED:
\vspace{-5pt}
\begin{equation}
PED=\frac{\Delta d/d}{\Delta \lambda/\lambda} \label{ped}
\vspace{-5pt}
\end{equation}
where $d$ denotes demand. To characterize PED more accurately, one needs to bear in mind that the overall PED should have the firm component of demand in it. Hence, in this setup, we set $d=d^{f}+v_{avg}(\lambda)$, where $d^{f}$ denotes the firm component of demand. We can observe in Figure \ref{figa} that the average PED is almost zero for prices that are considerably larger or smaller than the mean price, and only in the mid-portion of the plot (i.e. around the mean price) we notice a substantial average PED. One can verify using equation (\ref{ped}) with $d=d^{f}+v_{avg}(\lambda)$ that the average PED (i.e. the PED computed using the average storage response) depends on how the average storage response compares with the firm demand. Table \ref{tblped} shows the average PED around the mean price, using different values for $d^f$. %
\begin{table}[tp]%
\caption{Average price elasticity of demand around the mean price using different values for $d^f$}
\label{tblped}\centering%
\begin{tabular}{cc}
\hline\hline
 $d^f$  & Average Price Elasticity	\\ \hline 
 $\overline{v} $& -3.6		\\
	$3\overline{v}$ & -1.2	\\
	$ 8\overline{v}$ & -0.45		\\ \hline 
%With Penalties &  $ \overline{v}$ &-2.9		\\
%& $3\overline{v}$ &-0.96	\\
%&	$8\overline{v}$ &-0.36	\\ \hline \hline  
\end{tabular}
\vspace{-6pt}
\end{table}

\subsection{Price Responsiveness of a Storage System} \label{infelast}
\begin{figure}
[ptbh]
\begin{center}
\vspace{-17pt}
\includegraphics[scale=0.37]{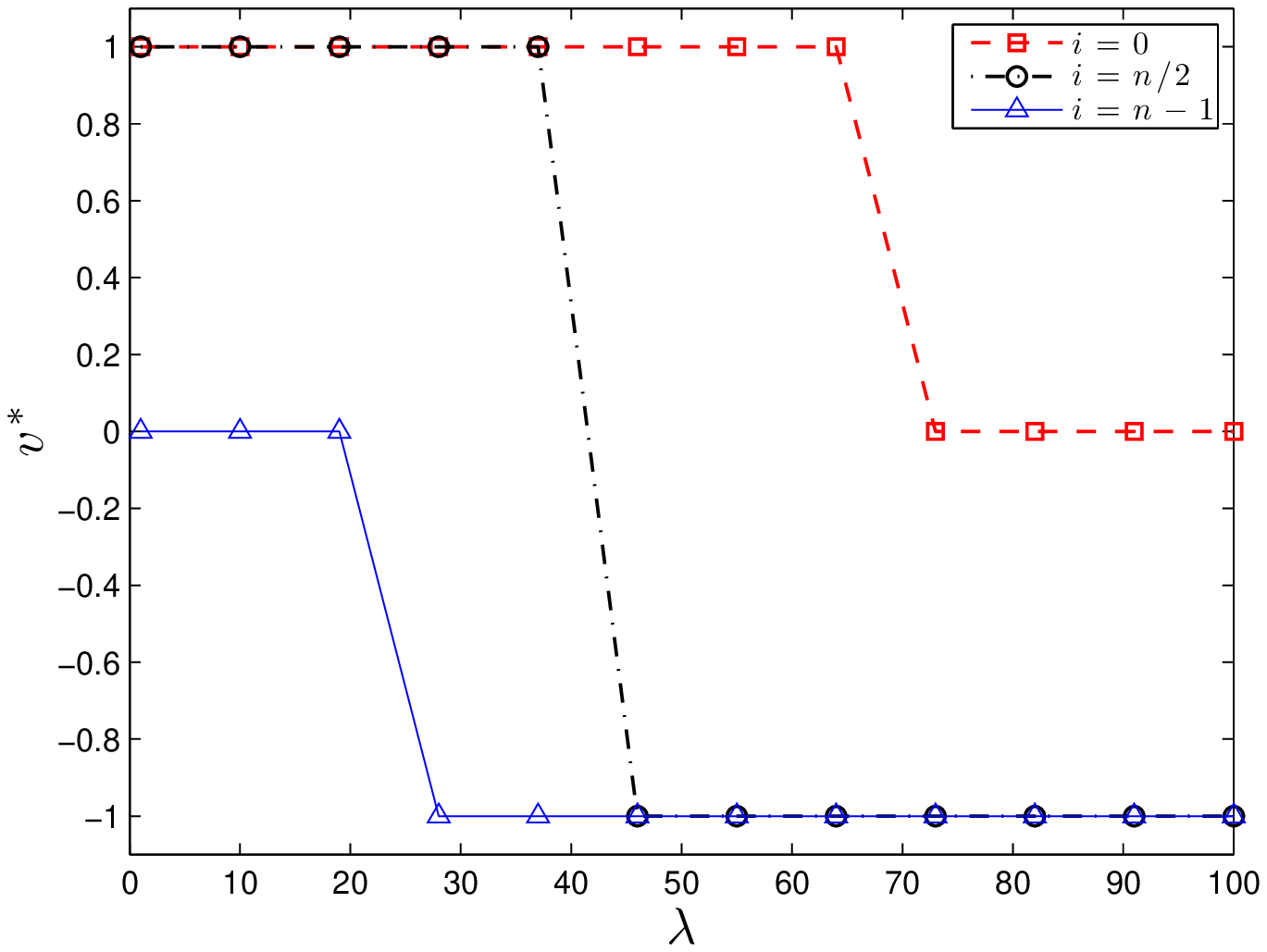}
\includegraphics[scale=0.37]{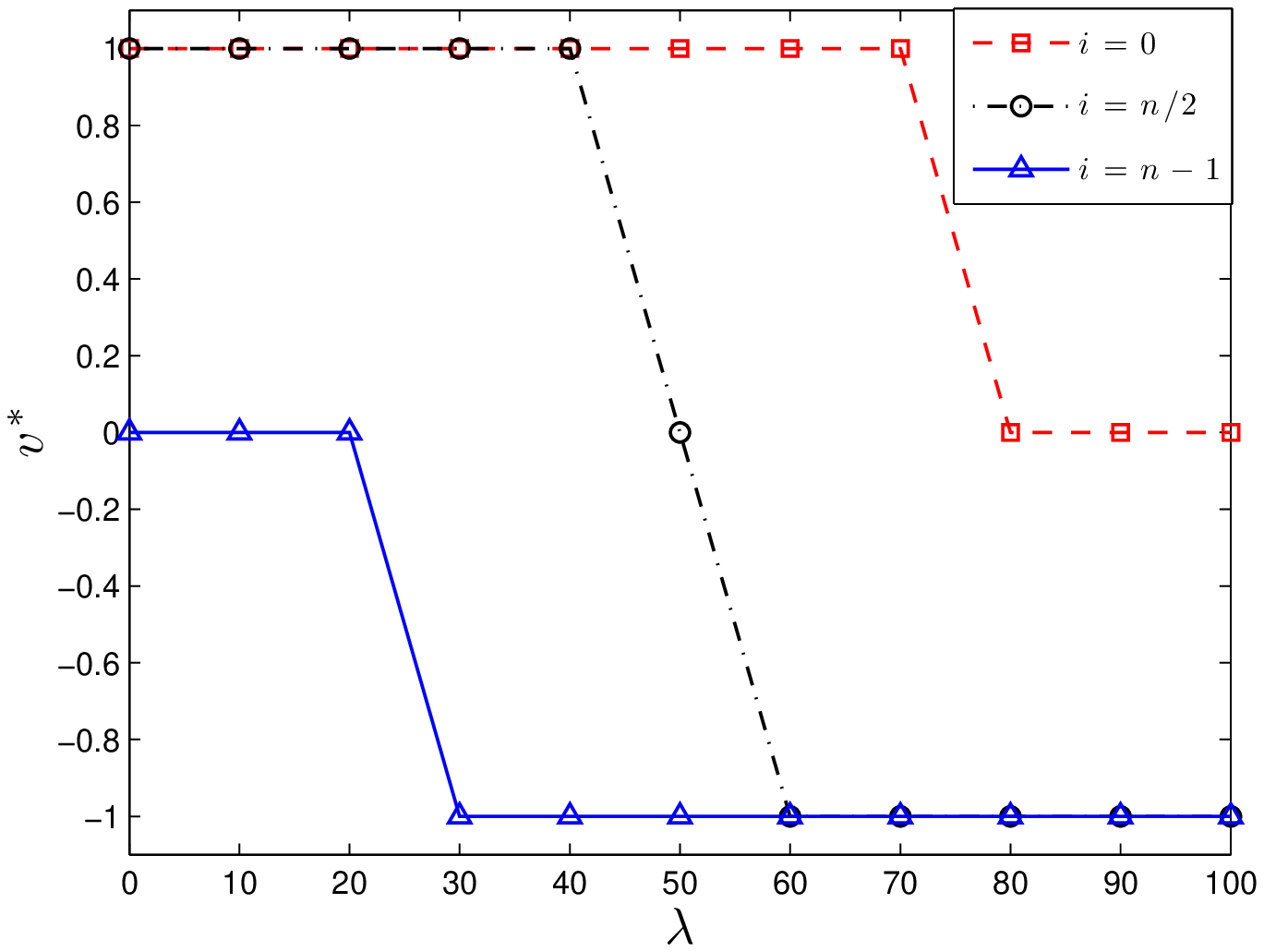}
\vspace{-3pt}
\caption{Storage response vs. price for three sample states for discretized truncated log-normal (left) and discretized uniform (right) price distributions, both with mean 50}
\label{figinf}
\end{center}
\vspace{-8pt}
\end{figure}

In this subsection, we present a computational framework for understanding the behavior of storage as a function of price and the amount of stored energy, and for characterization of the buy/sell phase transition region in the price-state plane. In order to eliminate stage dependence, we will consider the infinite-horizon version of the storage problem (\ref{storage}) and perform policy iteration (see, e.g.,  \citealt{Berts}) to numerically obtain a stationary (stage-independent) policy for purchase/sale as a function of both state and price. This section provides a qualitative picture of the structural characteristics of the behavior of storage, and a framework for estimating the PED as a function of the state. Herein, we will assume that the consumer starts with an empty storage, implying that the states would only take on integer multiples of the ramp constraint. We also use the same price distribution for all stages (i.e. we assume that prices are i.i.d). In our computations, we first use a discretized truncated log-normal price distribution, and then compare the results against the case of a discretized uniform distribution. For both distributions, we use a mean of about $\overline{\lambda}=50$ and a standard deviation of about $\sigma=30$, and also the same support. We set $\overline{v}=1$ and $n=10$. Figure \ref{figinf} illustrates how the storage response varies with price for three cases of the state (when the storage is empty ($i=0$), when the storage is half full ($i=n/2$), and when the storage is nearly full ($i=n-1$)).

Note in Figure \ref{figinf} the considerable effect of the storage state on storage response, compared to the small effect of the price distribution. More specifically, for both distributions, when the storage is empty, the optimal policy recommends purchasing from the grid even when the prices are somewhat above the mean price. Note that for the log-normal case, this policy change occurs at a slightly lower price because of the left skewness of the log-normal distribution.  
Though, when the storage is half full, we switch from the ``buy-it-all" policy to the ``sell-it-all policy" right at the mean price for the uniform distribution; for the log-normal distribution, this policy change occurs slightly before the mean price, which is again due to the left skewness of the log-normal distribution. Finally, when the storage is nearly full, for both distributions the optimal policy is to sell as much energy as possible for most prices, and to do nothing for the low prices. 

\begin{figure}
[ptbh]
\begin{center}
\vspace{-12pt}
\includegraphics[scale=0.37]{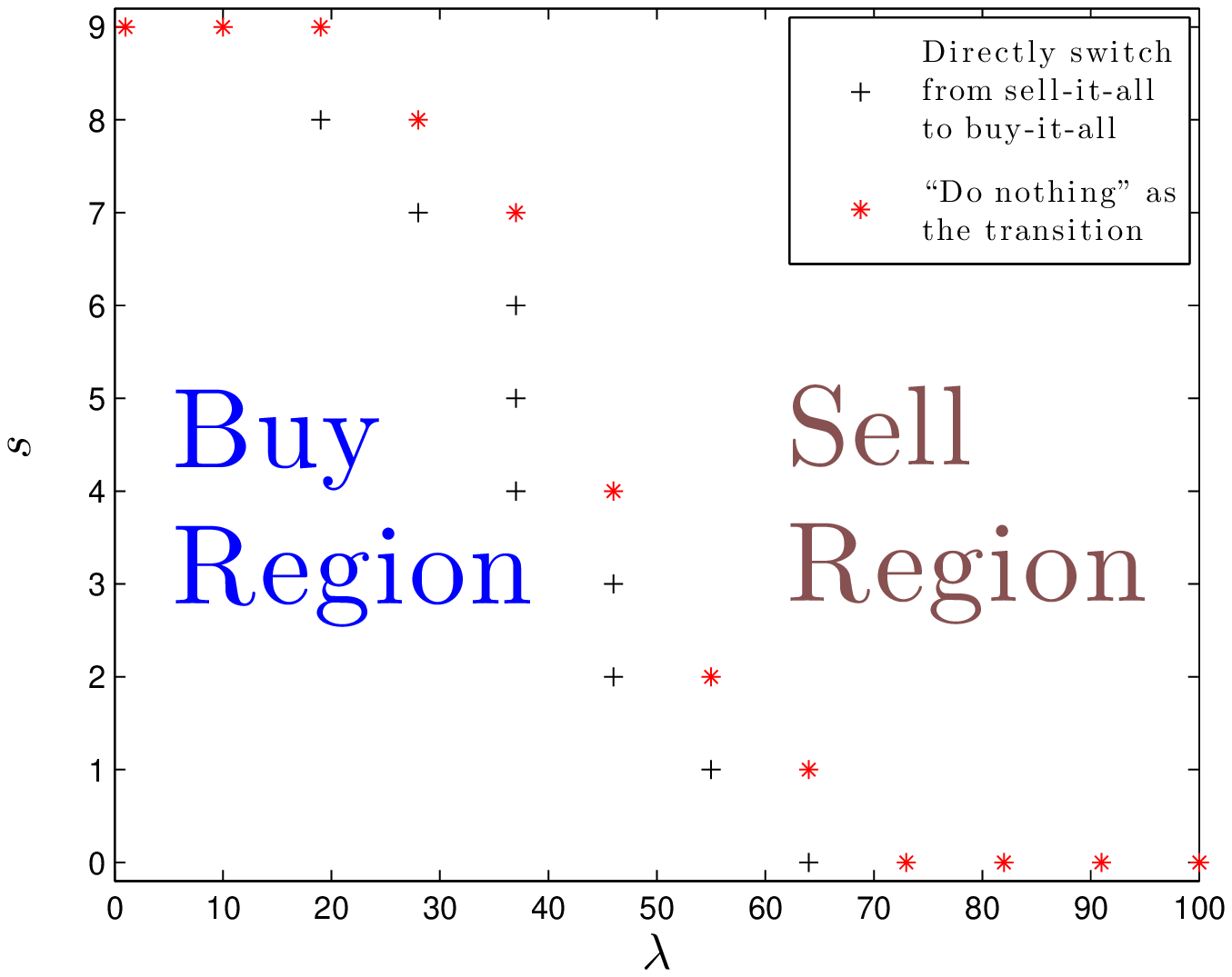}
\includegraphics[scale=0.385]{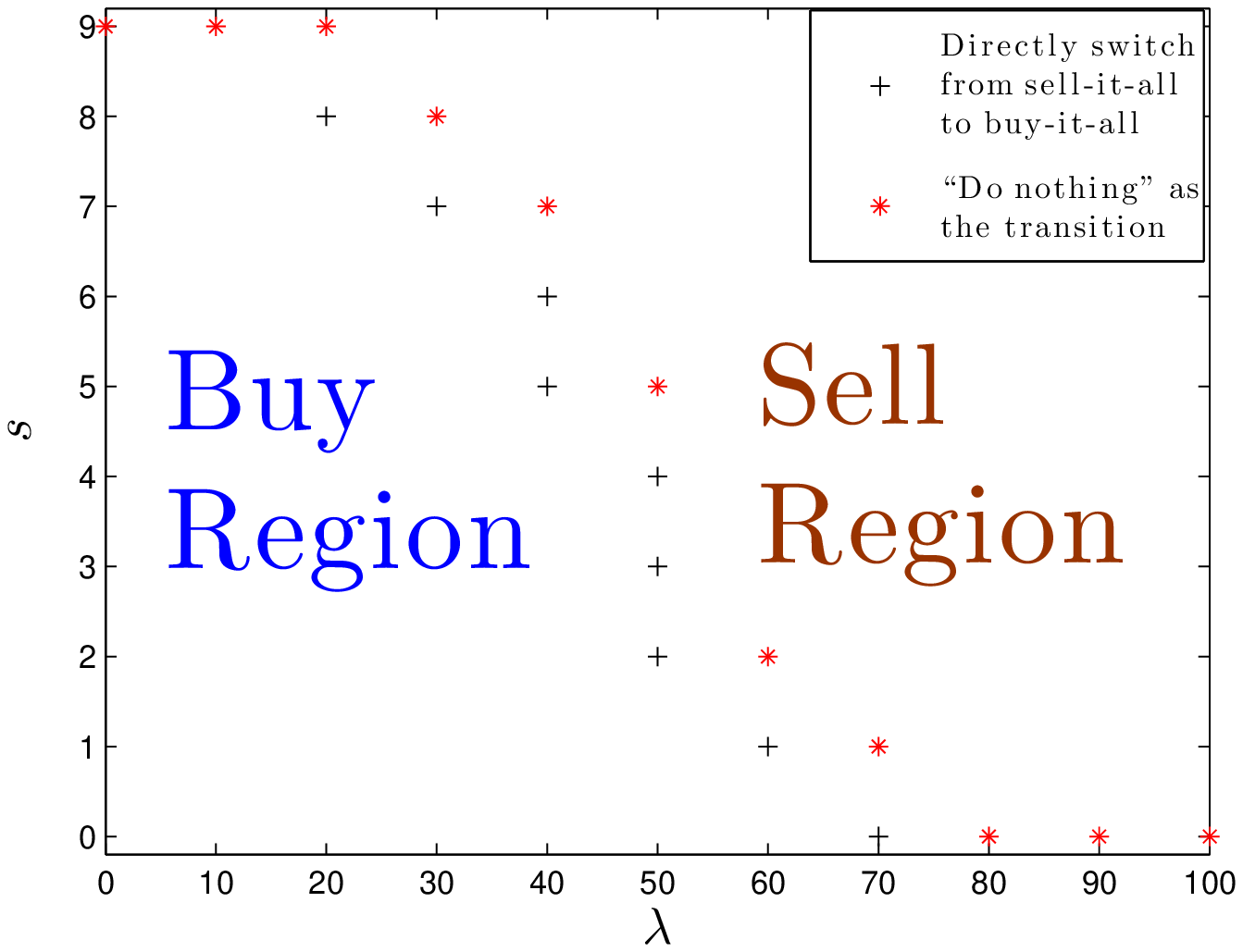}
%\vspace{-10pt}
\caption{Occurrence points of policy change from sell-it-all to buy-it-all for discretized truncated log-normal (left) and discretized uniform (right) price distributions}
\label{figinf2}
\end{center}
\vspace{-10pt}
\end{figure}

The transition points in the infinite-horizon policy from sell-it-all to buy-it-all on the $s-\lambda$ plane are shown in Figure \ref{figinf2}. Any point to the left of and/or below the transition points is a buying policy, which corresponds to $v^\ast(s,\lambda)=\overline{v}$, and any point to the right of and/or above the transition points is a selling policy, which corresponds to $v^\ast(s,\lambda)=-\overline{v}$. The plus signs show a direct transition from buying to selling when moving along the vertical axis, i.e., as storage state varies, unless they are immediately followed by a star on their right. The stars denote a transition through a ``Do Nothing" policy when moving along the horizontal axis, i.e., as price varies. Therefore, at the prices denoted by \textcolor{red}{*} we have $v^\ast=0$.
%\vspace{-6pt}
Figure \ref{figinf2} clearly illustrates the interplay between the state-dependence and the price-dependence of the storage response. Note also in Figure \ref{figinf2} the small effect of the price distribution on the storage response. The optimal policy for the log-normal case is slightly more shifted to left compared to that of the uniform case, which is again caused by the left skewness of the log-normal distribution.

We can now use the above results to characterize the PED as a function of the storage state. In order to compute the overall PED using  (\ref{ped}), once again we quantify the demand ($d$) such that it includes the firm component of the demand as well: 
\vspace{-5pt}
\begin{equation}
d=d^{f}+v^{\ast}(s,\lambda).\label{ovdem}
\vspace{-5pt}
\end{equation}
We have PED$(s)=0$ for all the points in the ``Buy Region" and in the ``Sell Region" because the storage response is constant in those regions. However, around the transition curve, the PED is non-zero because $\Delta d<0$ in that region. But $\Delta d$ can only take two values: $\Delta d=-2\overline{v}$ when there is a direct transition from the ``Buy Region" to the ``Sell Region", and $\Delta d=-\overline{v}$ when the transition is through a ``Do Nothing" policy. Hence, we have
\[
PED(s)=-\frac{2\overline{v}\lambda}{(d^{f}+v^\ast(s,\lambda))\Delta\lambda}~~\text{and}~~PED(s)=-\frac{\overline{v}\lambda}{(d^{f}+v^\ast(s,\lambda))\Delta\lambda}
\]
%\[
%P
%\]
around the points denoted by + and \textcolor{red}{*}, respectively, where, depending on the point at which we choose to compute PED, $v^{\ast}$ takes on one of the values in $\{-\overline{v},0,\overline{v}\}$. 

%%%%%%%%%%%%%%%%%%%%%%%%%%%%%%%%%%%%%%%%%%%%%%%%%%%%%%%%%%%%%%%%%%%%%%%%%%%%%%%%

\subsection{Impact of Market-Based Operation of Storage on the Required Reserves} \label{sec:reserves}
\begin{figure}
[ptbh]
\begin{center}
\vspace{-1.3in}
\includegraphics[scale=0.6]{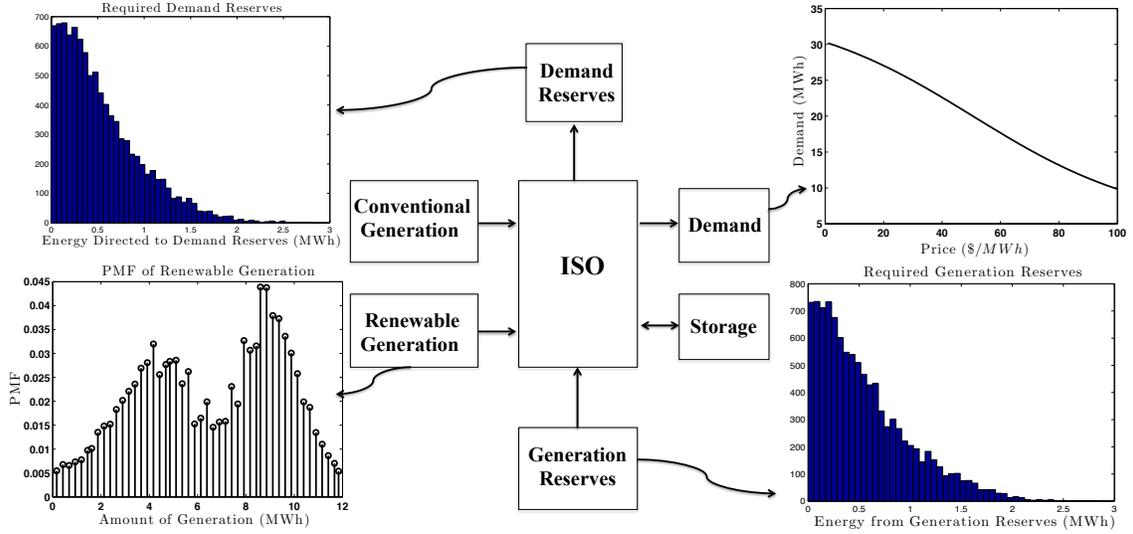}%{
\vspace{-1.2in}
\caption{Matching demand and supply by the ISO}
\label{iso_reser_chrt}
\end{center}
%\vspace{-15pt}
\end{figure}

In this subsection, we will evaluate how the optimal response from storage affects the amount of reserves (both generation and demand) required to guarantee that supply and demand can be matched. We assume that the Independent System Operator (ISO) primarily uses renewable generation with zero marginal cost and a conventional generation source with a quadratic cost to meet the demand. It is assumed that the overall demand consists of the storage response and a deterministic demand with both elastic and inelastic components. In our simulated setup, at the beginning of each pricing period, the ISO predicts the amount of renewable generation available during that period, possibly with some error. We assume that the amount of renewable generation in each period is i.i.d and is sampled from the bimodal distribution shown in Figure \ref{iso_reser_chrt}, which can for instance, correspond to a system with two renewable sources, one with high capacity and one with low capacity. In the event that during any time period, there is a shortfall in renewable generation compared to what was predicted, the ISO extracts from the generation reserves. Similarly, if there is an excess of renewable generation, the ISO would direct the excess generation to the demand reserves. We assume that the consumers can learn only the stationary distribution of prices (under the assumption of i.i.d prices), and not the exact mechanism by which prices are generated, and that the consumers are rational, and hence, they have no incentive to manipulate their storage response. Thus, the storage management policy discussed in Section \ref{infelast} is deemed optimal by the consumers. The setup is shown in Figure \ref{iso_reser_chrt}.  %(nrline{s}=5$).  
In this setup, at the beginning of each time period, the ISO determines the clearing price by solving the following equation for $\lambda_k$:
\vspace{-5pt}
\begin{equation}
a\lambda_{k}+\hat{w}_{k}= v^{\ast}(s,\lambda_{k})+d_{k}(\lambda_{k})
\vspace{-5pt}
\end{equation}
where $\hat{w}_{k}$ is the predicted amount of renewable generation during the $k$-th time period, $v^{\ast}(s,\lambda_{k})$ is the optimal response of storage obtained using policy iteration as discussed in Section \ref{infelast}, $d_{k}(\lambda_{k})$ is the deterministic demand, and $a\lambda_{k}$ is the optimal response of a supplier with quadratic cost $c(x)=\frac{1}{2a}x^2$, to a given price $\lambda_k$. For the purpose of numerical computations, we assume that $d_{k}(\lambda_{k})$ is a logistic function (as shown in Figure \ref{iso_reser_chrt})
 \begin{equation}
 d_{k}(\lambda_{k})=5~\text{MWh}+\frac{(30~\text{MWh})e^{-\lambda_{k}/30}}{e^{-\lambda_{k}/30}+e^{-\overline{\lambda}/30}}.
 \end{equation}
The choice of a logistic function for characterizing the deterministic demand was due to the fact that qualitatively similar functions have been used in the past for modeling price responsive demand in electricity markets (see, e.g., \citealt{Carr}), this choice is otherwise arbitrary and not central to our analysis or conclusions. The choice of the numerical values of the parameters of the logistic function is reasonable for the hourly load in a small power system. For all $k$ we set $\overline{s}/\overline{v}=n=5$, the average price $\overline{\lambda}=50~\$/\text{MWh}$, and vary the ramp constraint for experimentation, looking at three cases: $\overline{v}=0.25 ~\text{MW}$, $\overline{v}=0.5~ \text{MW}$, and $\overline{v}=1~ \text{MW}$ such that the storage capacity $\overline{s}$ is roughly $4.17\%,~8.33\%$, and $16.67\%$, respectively, of the maximum possible demand in the absence of storage. We also set $a=4~(\text{MWh})^{2}/\$$ so that the average price in the simulated setup is about 50 $\$\text{/MWh}$, which is close to the typical average hourly prices in various electricity markets. Note that since we are fixing $n$, all the storage devices across the grid become synchronized in the steady state, and hence, we can use one large value for $\overline{v}$ which serves as a proxy for the cumulative ramp rate constraint across the grid.

We further assume that the prediction of the amount of renewable generation available in each period has an error that is sampled from a truncated normal distribution with zero mean and standard deviation 0.5 MWh, and that the ISO estimates the storage state with an error that has a discrete uniform distribution over the support $\{-\overline{v},0,\overline{v}\}$. We compute the percentage change in the amount of reserves needed in the presence of storage compared to the case of no storage, for each storage capacity mentioned above. In each scenario, we first examine the 100$\%$ reliability level (i.e. when demand and supply are guaranteed to match roughly $100\%$ of the time) and then we repeat the computations for $99\%$ and $98\%$ reliability levels.
\begin{table}[tp]%
\caption{Change in the required reserves for various values of ramp and reliability levels}
\label{reserve_res}\centering%
\begin{tabular}{cccc}
\hline\hline
$\overline{s}$/max hourly dem  &  Reliability & Change in demand reserves &Change in generation reserves	\\ \hline 
\large$\frac{1}{24}$\normalsize	 & $100\%$&$7.0\%\uparrow$&$9.9\%\uparrow$		\\
& $99\%$& 3.1\%$\uparrow$&$2.7\%\uparrow$		\\
 & $98\%$& 2.4\%$\uparrow$&3.8\%$\uparrow$	\\ \hline 
\large$\frac{1}{12}$\normalsize & $100\%$&22\%$\uparrow$&$24\%\uparrow$		\\
& $99\%$&13$\%\uparrow$&$15\%\uparrow$		\\
 & $98\%$&$12\%\uparrow$&$14\%\uparrow$	\\ \hline 
\large$\frac{1}{6}$\normalsize 	 &$100\%$&40$\%\uparrow$&$42\%\uparrow$		\\
& $99\%$&$43\%\uparrow$ &$45\%\uparrow$		\\
 & $98\%$&$45\%\uparrow$&$46\%\uparrow$	\\ \hline 
\end{tabular}
\end{table}

As the simulation results in Table \ref{reserve_res} suggest, market based operation of storage may require an increase in the amount of required reserves (both demand and supply reserves). The larger the amount of energy that can be stored, or extracted from storage, the higher the amount of demand and supply reserves needed. Note that in the largest case of storage reported in Table \ref{reserve_res}, the storage capacity is only one sixth of the maximum demand, while the reserves need to be expanded by about 40$\%$ to accommodate the integration of this much storage (the histograms of the required reserves for this case of storage capacity are shown in Figure \ref{iso_reser_chrt}). 

This effect is mainly caused by the threshold-based, state dependent response of storage as discussed in Section \ref{infelast} (Figure \ref{figinf2}). When the prices are high due to lack of renewable generation, the storage system does not necessarily respond to the price incentive the way the ISO would want it to, because even if the prices are above the average, the optimal policy for the storage could be to buy if the storage level is low. Hence, in the event that the actual amount of renewable generation falls short of what was predicted, if the ISO also makes an error in estimating the storage state, a considerable amount of energy may be extracted from the grid by the storage device, which would force the ISO to use the generation reserves even more. Analogously, if there is an excess of renewable generation, the ISO may be forced to use the demand reserves even more if the storage device feeds energy back into the grid. A higher penetration of storage, and hence, a higher aggregate ramp rate amplifies this issue because the impact of erroneous state estimations would be even more profound. However, the generation reserves are typically fast generators with high economic and environmental costs. Hence, the ISO may need to regulate market-based operation of storage devices to prevent such impacts, and/or design mechanisms for pricing these externalities. Also, the ISO may need to design mechanisms that guarantee its access to the exact value of the storage state, because the intense interplay between state and price is such that estimating the storage state even with a small error can heavily impact the storage response.

%%%%%%%%%%%%%%%%%%%%%%%%%%%%%%%%%%%%%%%%%%%%%%%
%%%%%%%%%%%%%%%%%%%%%%%%%%%%%%%%%%%%%%%%%%%%%%%%%

%%%%%%%%%%%%%%%%%%%%%%%%%%%%%%%
%%%%%%%%%%%%%%%%%%%%%%%%%%%%%%%%%%%%%%%%%%%%%%%%%%%%%%%%%%%%%%%%%%%%%%%%%%%%%%%%
%\vspace{-1pt}
\section{Conclusion}
%\vspace{-2pt}
%%%%%%%%%%%%%%%%%%%%%%%%%%%%%%%%%%%%%%%%%%%%%%%%%%%%%%%%%%%%%%%%%%%%%%%%%%%%%%%%
%%%%%%%%%%%%%%%%%%%%%%%%%%%%%%%%%%%%%%%%%%%%%%%%%%%%%%%%%%%%%%%%%%%%%%%%%%%%%%%%
%%%%%%%%%%%%%%%%%%%%%%%%%%%%%%%%%%%%%%%%%%%%%%%%%%%%%%%%%%%%%%%%%%%%%%%%%%%%%%%%
\label{sec:conclusions}
In this paper, we proposed a dynamic model for optimal control of storage under ramp constraints and exogenous, stochastic prices. We derived the associated optimal policy and value function, and gave explicit formulas for their computation. Moreover, we derived an analytical upperbound on the long-term average economic value of storage, which is valid for any achievable realization of prices over a fixed support, and highlights the dependence of the value on ramp constraints and capacity. This result can be useful in assessing viability of investment in electricity storage. We also showed that while the value of storage is a non-decreasing function of price volatility, due to finite ramping rates, the value of storage saturates quickly as capacity increases, regardless of price volatility. We highlighted the dependence of the response of storage to prices on the internal state of storage, and also, by averaging out state and stage dependence, we showed that in expectation, storage induces a considerable amount of price
elasticity near the average price. We also showed that if the ISO does not have perfect information about the exact value of the storage state, the reserves may need to be expanded to accommodate market-based operation of storage.

Our results provide insight into learning the behavior of storage, particularly modeling and estimating the response of a ramp-constrained storage system when used as an arbitrage mechanism. We used price data from real-time wholesale markets to examine the sensitivity of the optimal policy and the value of storage to our assumption of independent prices. The relatively high competitive ratios that we found when testing the optimal policy suggest that the sensitivity of the policy to the assumption of independent prices is low. The competitive ratios found in this paper by applying the optimal policy can be used as a benchmark for comparison with other results that assume more complicated models of the stochastic price process.  %

We did not directly consider the effects of inefficiencies such as conversion losses and aging on the optimal policy. While such inefficiencies certainly limit the economic value of storage, an important question to ask is what is their effect on the response of storage to prices? For the case of battery storage, aging is proportional to the amount of current withdrawn or injected. Intuitively, this would make buying energy for, or selling energy from storage less profitable at moderate prices. Within the class of threshold policies, the corresponding optimal selling thresholds would be higher and the buying thresholds would be lower than what we have derived in this paper and the ``do nothing" range would be wider. Qualitatively, this would mean an overall choppier response from storage, with high elasticity over a narrow range of prices and low elasticity over a wide range, which would be undesirable from a system operation and reliability point of view. 

While our paper largely focused on the economic value of storage, it is important to recognize and quantify the environmental and the reliability value of storage. With proper control policies, storage can help matching stochastic supply with demand, improving system frequency and voltage profiles, and possibly mitigating large blackouts. The development of a systematic framework for quantifying the value of storage, the trade-offs between reliability, environmental, and economic value of storage, and the design of real-time pricing and market mechanisms for optimally striking these trade-offs are fundamental and important directions for future research.

\bibliographystyle{nonumber}

%\vspace{-2pt}
\newpage
\section{Appendix}
{\bf Proof of Theorem \ref{storageTHM}:}
\proof{Proof.}
This proof proceeds by induction. Let us for the moment assume that the value function $V_{k}(\cdot)\overset{\text{def}}= \mathbf{E}\left[  J_{k}\left(  \cdot \right)  \right]  $ has the form defined in (\ref{valuefunc}).\\
From the dynamic programming algorithm, for $k<N$,
we have
\begin{equation*}
J_{k}\left(  s_{k}\right)  =h_{k}\left(  s_{k}\right)+\min_{v_{k}\in \lbrack \max(-s_{k},-\overline{v}),\overline{v}] }\left\{  \lambda_{k}v_{k}
+\mathbf{E}\left[  J_{k+1}\left(  s_{k}+v_{k}\right)  \right]  \right\}, \tag{A.1}
\label{Nm1stepStor}
\end{equation*}
where the penalty functions $h_{k}(s_{k})$ are as defined in (\ref{penal}).\\
Then, we use the general form in (\ref{valuefunc}) for $\mathbf{E}\left[  J_{k}\left(  s_{k}\right)  \right] $ and the state evolution equation in (\ref{thestates}), and let $i\in\overline{\mathbb{Z}}_{+}$ be such that $v_{k}+s_{k}\in\left[i\overline
{v},(i+1)\overline{v}\right)$. Applying the induction step to (\ref{Nm1stepStor}) we obtain:
\begin{equation*}
J_{k}(s_{k})=h_{k}\left(  s_{k}\right)+\min_{v_{k}\in\left[  \max( -s_{k},-  \overline{v}),   \overline{v}\right]  }~\{  \lambda_{k}v_{k}
- t^{i}_{k+1}(s_{k}+v_{k})+e^{i}_{k+1}\}  \tag{A.2}\label{Nm2stepStor}
\end{equation*}
Solving the optimization problem in (\ref{Nm2stepStor}) yields the optimal policy shown in (\ref{optpol}).
%%\]
Now we can plug in the optimal policy to obtain $J_{k}\left(  s_{k}\right)$. Let $i\in\overline{\mathbb{Z}}_{+}$ be such that $s_{k}\in\left[i\overline
{v},(i+1)\overline{v}\right)$. If $0\leq s_{k} <  \overline{v}$, then
\[
J_{k}\left(  s_{k}\right)=\left\{
\begin{array}
[c]{lcl}%
\!\!(h^{0}_{k}-\lambda_{k})s_{k}+e^{0}_{k+1}+c^{0}_{k} & \text{if}& t^{0}_{k+1}<\lambda_{k} \vspace{0.05in}\\
\!\!(\lambda_{k}-t^{1}_{k+1})  \overline{v}%
+(h^{0}_{k}- \lambda_{k})s_{k}  +e^{1}_{k+1}+c^{0}_{k}& \text{if}&t^{1}_{k+1}<\lambda_{k} \leq t^{0}_{k+1} \vspace{0.05in}\\
\!\!(\lambda_{k}-t^{1}_{k+1})  \overline{v}%
+(h^{0}_{k}- t^{1}_{k+1})s_{k}  +e^{1}_{k+1}+c^{0}_{k}& \text{if}&\lambda_{k}~~\leq t^{1}_{k+1}\vspace{0.05in}\\
\end{array}
\right.
\]%
and if $s_{k}\geq \overline{v}$, then
\[
J_{k}\left(  s_{k}\right)\!=\!\!\left\{
\begin{array}
[c]{lcl}
\!\!-(\lambda_{k}- t^{i-1}_{k+1})  \overline{v}%
+(h^{i}_{k}- t^{i-1}_{k+1})s_{k}  +e^{i-1}_{k+1}+c^{i}_{k}& \text{if}& t^{i-1}_{k+1}<\lambda_{k} \vspace{0.05in}\\
\!\!(\lambda_{k}-t^{i}_{k+1})i  \overline{v}%
+(h^{i}_{k}- \lambda_{k})s_{k}  +e^{i}_{k+1}+c^{i}_{k}& \text{if}&t^{i}_{k+1}<\lambda_{k}\leq t^{i-1}_{k+1} \vspace{0.05in}\\
\!\!(\lambda_{k}-t^{i+1}_{k+1})(i+1)  \overline{v}%
+(h^{i}_{k}- \lambda_{k})s_{k}  +e^{i+1}_{k+1}+c^{i}_{k}& \text{if}&t^{i+1}_{k+1}<\lambda_{k}\leq t^{i}_{k+1} \vspace{0.05in}\\
\!\!(\lambda_{k}- t^{i+1}_{k+1}) \overline{v}%
+(h^{i}_{k}- t^{i+1}_{k+1})s_{k} +e^{i+1}_{k+1}+c^{i}_{k} & \text{if}&\lambda_{k}~~\leq t^{i+1}_{k+1}
\end{array}
\right.
\]%

\noindent Let us recall that for a piecewise function $f$ defined according to
\[
f\left(  x\right)  =\left\{
\begin{array}
[c]{ccc}%
\!\!f_{1}\left(  r\right)  x+g_{1}\left(  r\right)  & \text{if} & r\leq \alpha \vspace{0.05in}\\
\!\!f_{2}\left(  r\right)  x+g_{2}\left(  r\right)  & \text{if} & r>\alpha
\end{array}
\right.
\]
\vspace{-5pt}
we have
\begin{align}
\vspace{-11pt}
\mathbf{E}\left[  f\left(  x\right)  \right]   &  =\mathbf{E}\left[  f\left( 
x\right)  |r\leq\alpha\right]  \mathbf{P}\left(  r\leq\alpha\right)  +\mathbf{E}%
\left[  f\left(  x\right)  |r>\alpha\right]  \mathbf{P}\left(  r>\alpha\right)\nonumber
\\
&  ~~~+\mathbf{E}\left[  g_{1}\left(  r\right)  |r\leq \alpha\right]  \mathbf{P}%
\left(  r\leq\alpha\right)  +\mathbf{E}\left[  g_{2}\left(  r\right)
|r>\alpha\right]  \mathbf{P}\left(  r>\alpha\right) \tag{A.3} \label{totalExp}
\vspace{-11pt}
\end{align}
Also, note that:
\begin{equation}
\vspace{-5pt}
\mathbf{E}\left[  f\left(  r\right)  |r\leq\alpha\right]  \mathbf{P}\left(
r\leq \alpha\right)  ={\displaystyle\sum\limits_{\theta=r_{\min}}^{\alpha}}
f\left(  r\right)  P_{R}\left(  \theta\right)  \tag{A.4}\label{condExp}
\vspace{-5pt}
\end{equation}
\\
Now, let us apply (\ref{totalExp}) and (\ref{condExp}) to the equations derived above for $J_{k}\left(  s_{k}\right)$, and compute $ \mathbf{E}\left[  J_{k}\left(  s_{k}\right)  \right] $ for $k<N$, which leads to the following results:\\
\[
\begin{array}
[c]{lcl}
\text{if}~0\leq s_{k} <  \overline{v},~\text{then}\\
\mathbf{E}\left[  J_{k}\left(  s_{k}\right)  \right]=e^{0}_{k} -s_{k}[ t^{1}_{k+1} F_{k}\left(  t^{1}_{k+1}\right)
 -h^{0}_{k}+
{\displaystyle\sum\limits_{ t^{1}_{k+1}<\theta \leq \lambda^{\max}_{k}}^{}}
\theta P_{k}(  \theta)]  \\[0.05in]
\\
\text{if}~ s_{k}\geq \overline{v}, \text{then}\\
\mathbf{E}\left[  J_{k}\left(  s_{k}\right)  \right]=e^{i}_{k} -s_{k}[t^{i-1}_{k+1}\left(  1-F_{k}\left(  t^{i-1}_{k+1}\right)  \right)
+t^{i+1}_{k+1}F_{k}\left(  t^{i+1}_{k+1}\right)  -h^{i}_{k}+{\displaystyle\sum\limits_{t^{i+1}_{k+1}<\theta \leq t^{i-1}_{k+1}}^{}}
\theta P_{k}(  \theta)]  \\
\end{array}
\]
where $e^{i}_{k}$ denotes the sum of the terms that have not been multiplied by $s_{k}$. Hence, the thresholds for $k<N$ and $i\in\mathbb{Z}_{+}$ are given by: 
\[
\begin{array}
[c]{lcl}
t^{0}_{k}= t^{1}_{k+1} F_{k}\left(  t^{1}_{k+1}\right)
 -h^{0}_{k}+
{\displaystyle\sum\limits_{ t^{1}_{k+1}<\theta \leq \lambda^{\max}_{k}}^{}}

\theta P_{k}(  \theta)  \\[0.05in]
t^{i}_{k}=t^{i-1}_{k+1}\left(  1-F_{k}\left(  t^{i-1}_{k+1}\right)  \right)
+t^{i+1}_{k+1}F_{k}\left(  t^{i+1}_{k+1}\right)  -h^{i}_{k}\!+\!
{\displaystyle\sum\limits_{t^{i+1}_{k+1}<\theta \leq t^{i-1}_{k+1}}^{}}
\theta P_{k}(  \theta) \\
\end{array}
\]
After some algebra, the above results can be written in the form presented in (\ref{thresh}). The derivation of the intercepts $e^{i}_{k}$ follows the same procedure.

The next step is to verify, using induction, that the thresholds at each stage (i.e. $t^{i}_{k}$) are non-increasing functions of $i$. Considering that $t^{0}_{k-1}$ has a different general form than $t^{i}_{k-1}$ for $i>0$, we first need to show that $t^{1}_{k-1} \leq t^{0}_{k-1}$ assuming that $t^{i+1}_{k}\leq t^{i}_{k}$ for all $i$. Hence, we need to first show that\\ \fontsize{11.2}{12}%\smallsize
 $ t^{1}_{k}F_{k-1}\left(  t^{1}_{k}\right)
  -h^{0}_{k-1}+
{\displaystyle\sum\limits_{ t^{1}_{k}<\theta \leq \lambda^{\max}_{k-1}}^{}}
\theta P_{k-1}\left(  \theta\right)\geq t^{0}_{k}\left(1- F_{k-1}\left(  t^{0}_{k}\right)
\right) + t^{2}_{k} F_{k-1}\left(  t^{2}_{k}\right)
 - h^{1}_{k-1}+
{\displaystyle\sum\limits_{t^{2}_{k}<\theta \leq t^{0}_{k}}^{}}
\theta P_{k-1}\left(  \theta\right) $ \\ \\
\normalsize Knowing that $h^{1}_{k-1}\geq h^{0}_{k-1}$, we can remove $-h^{1}_{k-1}$ and $-h^{0}_{k-1}$ from both sides. Then, by writing the cumulative distribution functions as summations, and rearranging some terms, the above can be rewritten as:\\  \fontsize{11.2}{12}
 $t^{0}_{k}\!\!\!\!\!\! {\displaystyle\sum\limits_{\lambda^{\min}_{k-1}\leq \theta \leq t^{0}_{k}}^{}}
P_{k-1}\left(  \theta\right)
 +\!\!\!\!\!\! {\displaystyle\sum\limits_{\lambda^{\min}_{k-1}\leq \theta \leq t^{1}_{k}}^{}}
 P_{k-1}\left(  \theta\right)t^{1}_{k} +\!\!\!\!\!\!
 {\displaystyle\sum\limits_{ t^{1}_{k}<\theta \leq \lambda^{\max}_{k-1}}^{}}
\theta P_{k-1}\left(  \theta\right) \geq t^{0}_{k}+\!\!\!\!\!\! {\displaystyle\sum\limits_{\lambda^{\min}_{k-1}\leq \theta \leq t^{2}_{k}}^{}}
 P_{k-1}\left(  \theta\right) t^{2}_{k} 
  +\!\!\!\!\!\!{\displaystyle\sum\limits_{t^{2}_{k}<\theta \leq t^{0}_{k}}^{}}
\theta P_{k-1}\left(  \theta\right)$\\ \\
\normalsize By breaking each summation into disjoint intervals, and factoring all the terms in the same interval and merging them into one summation, the above can be rewritten as:\\ \fontsize{11.2}{12}
 $t^{0}_{k} ~\leq \!\!\!\!\!\!{\displaystyle\sum\limits_{\lambda^{\min}_{k-1}\leq \theta \leq t^{2}_{k}}^{}}
 P_{k-1}\left(  \theta\right)( t^{0}_{k}+t^{1}_{k}-t^{2}_{k}) 
  +
\!\!\!\!\!\!{\displaystyle\sum\limits_{t^{2}_{k}<\theta \leq t^{1}_{k}}^{}}
 P_{k-1}\left(  \theta\right)  (t^{0}_{k}+t^{1}_{k}-\theta)+\!\!\!\!\!\!{\displaystyle\sum\limits_{t^{1}_{k}< \theta \leq t^{0}_{k}}^{}}
P_{k-1}\left(  \theta\right)t^{0}_{k}
 +\!\!\!\!\!\! \displaystyle\sum\limits_{t^{0}_{k}< \theta \leq {\lambda^{\max}_{k-1}}^{}}
\theta P_{k-1}\left(  \theta\right)  $\\ \\
\normalsize We can see that in the equation on the right hand side (RHS) of the inequality shown above, all the terms that have been multiplied by $P_{k-1}\left(  \theta\right)$ inside the summations are greater than or equal to $t^{0}_{k}$; we also know that ${\displaystyle\sum\nolimits_{\theta= \lambda^{\min}_{k-1}}^{\lambda^{\max}_{k-1}}}P_{k-1}\left(  \theta\right) = 1 $. Hence, we verify by inspection that the RHS equation in the inequality shown above is always greater than or equal to $t^{0}_{k}$.

Finally, we have assigned a non-positive cost of $-\hat{t}$ to each unit of energy left in storage at stage $N$ (i.e. $\mathbf{E}\left[  J_{N}\left(  s_{N}\right)  \right]=-\hat{t}s_{N}$ for $s_{N} \leq \overline{s}$). Also, a very small (negative) value is assigned to the thresholds for $i\geq n$ at stage $N$, to make sure we will not exceed the storage capacity in this stage. Hence, $t^{i}_{N}$ is a non-increasing function of $i$, and $ t^{i+1}_{N} \leq t^{i}_{N}$ is satisfied. 

We must now verify that $t^{i+1}_{k-1}\leq t^{i}_{k-1}$ for $i\in\mathbb{Z}_{+}  $ assuming that $t^{i}_{k}$ is a non-increasing function of $i$. So, we need to verify:\\ \fontsize{11.2}{12}
$ t^{i}_{k}\left(1- F_{k-1}\left(  t^{i}_{k}\right)
\right) + t^{i+2}_{k}F_{k-1}\left(  t^{i+2}_{k}\right)
 -h^{i+1}_{k-1}+\!\!\!\!\!\!
{\displaystyle\sum\limits_{t^{i+2}_{k}<\theta \leq t^{i}_{k}}^{}}\!\!\!\!\!\!
\theta P_{k-1}\left(  \theta\right) ~\leq~t^{i-1}_{k}\left(1- F_{k-1}\left(  t^{i-1}_{k}\right)
\right) \\~~~~~~~~~~~~~~~~~~~~~~~~~~~~~~~~~~~~~~~~~~~~~~~~~~~~~~~~~~~~~~~~~~~~~~~~~~~~~~~+ t^{i+1}_{k} F_{k-1}\left(  t^{i+1}_{k}\right)+\!\!\!\!\!\!{\displaystyle\sum\limits_{t^{i+1}_{k}<\theta \leq t^{i-1}_{k}}^{}}\!\!\!\!\!\!
\theta P_{k-1}\left(  \theta\right)  -h^{i}_{k-1}$\\ \\
\normalsize Knowing that $h^{i+1}_{k-1}\geq h^{i}_{k-1}$, we can remove $-h^{i+1}_{k-1}$ and $-h^{i}_{k-1}$ from both sides. Then, by writing the cumulative distribution functions as summations, and rearranging some terms, the above can be rewritten as:\\  \fontsize{11.2}{12}
$ t^{i}_{k}+ t^{i-1}_{k}\!\!\!\!{\displaystyle\sum\limits_{ \lambda^{\min}_{k-1}\leq \theta \leq t^{i-1}_{k}}^{}}
\!\!\!\!P_{k-1}\left(  \theta\right)+ t^{i+2}_{k}\!\!\!\!{\displaystyle\sum\limits_{ \lambda^{\min}_{k-1}\leq \theta \leq t^{i+2}_{k}}^{}}
\!\!\!\!P_{k-1}\left(  \theta\right)+
\!\!\!\!{\displaystyle\sum\limits_{t^{i+2}_{k}<\theta \leq t^{i}_{k}}^{}}
\theta P_{k-1}\left(  \theta\right)  ~\leq~t^{i-1}_{k}+ t^{i}_{k}{\displaystyle\sum\limits_{ \lambda^{\min}_{k-1}\leq \theta \leq t^{i}_{k}}^{}}
P_{k-1}\left(  \theta\right)$\\$~~~~~~~~~~~~~~~~~~~~~~~~~~~~~~~~~~~~~~~~~~~~~~~~~~~~~~~~~~~~~~~~~~~~~~~~~~~~~~~~~~~~~~~~~~~~~~~~~ +t^{i+1}_{k}\!\!\!\!\!\!{\displaystyle\sum\limits_{ \lambda^{\min}_{k-1}\leq \theta \leq t^{i+1}_{k}}^{}}
\!\!\!\!P_{k-1}\left(  \theta\right) +\!\!\!\!{\displaystyle\sum\limits_{t^{i+1}_{k}<\theta \leq t^{i-1}_{k}}^{}}
\!\!\!\!\theta P_{k-1}\left(  \theta\right)  $\\ 
%\th
\newline \normalsize By taking all the summations to the right hand side and taking $t^{i}_{k}$ to the left hand side of the inequality, breaking each summation into disjoint intervals, and factoring all the terms in the same interval and merging them into one summation, the above can be rewritten as:\\  \fontsize{11.2}{12}
$t^{i-1}_{k}-t^{i}_{k} \geq {\displaystyle\sum\limits_{ \lambda^{\min}_{k-1}\leq \theta \leq t^{i+2}_{k}}^{}}
P_{k-1}\left(  \theta\right)(t^{i-1}_{k}-t^{i}_{k}-(t^{i+1}_{k}-t^{i+2}_{k}))+ {\displaystyle\sum\limits_{t^{i+2}_{k}<\theta \leq t^{i+1}_{k}}^{}}
P_{k-1}\left(  \theta\right)( t^{i-1}_{k}-t^{i}_{k}-(t^{i+1}_{k}-\theta))$\\$~~~~~~~~~~~~~~~~~+\!\!\!\!\!\!{\displaystyle\sum\limits_{t^{i+1}_{k}<\theta \leq t^{i}_{k}}^{}}
P_{k-1}\left(  \theta\right)(t^{i-1}_{k}-t^{i}_{k})+ {\displaystyle\sum\limits_{t^{i}_{k}<\theta \leq t^{i-1}_{k}}^{}}
P_{k-1}\left(  \theta\right)(t^{i-1}_{k}- \theta ) $\\  \\
\normalsize We can see that in the equation on the RHS of the inequality shown above, all the terms that have been multiplied by $P_{k-1}\left(  \theta\right)$ inside the summations are less than or equal to $t^{i-1}_{k}-t^{i}_{k}$; we also know that ${\displaystyle\sum\nolimits_{\theta= \lambda^{\min}_{k-1}}^{t^{i-1}_{k}}}P_{k-1}\left(  \theta\right) \leq 1 $. So, given that the summations in the RHS do not overlap, we can verify by inspection that the RHS equation in the inequality shown above is always less than or equal to $t^{i-1}_{k}-t^{i}_{k}$.

The last step is to show that the value function is a continuous function. This can be done by induction. We have defined $\mathbf{E}\left[  J_{N}\left(  s_{N}\right)  \right]$ in such a way that it is convex. We also have defined $h_{k}(s_k)$ to be convex for all $k$. Looking at the equations of $J_{k}(s_{k})$ in (\ref{Nm1stepStor}) and (\ref{Nm2stepStor}), given that $\mathbf{E}\left[  J_{k+1}\left(  s_{k+1}\right)  \right]$ is convex, one can observe that the continuity of  $J_{k}(s_{k})$ is satisfied because $v_{k}^\ast$ is a continuous function of $s_k$. Similarly, considering the equations obtained for $v_{k}^\ast$, the continuity of $\mathbf{E}\left[  J_{k}\left(  s_{k}\right)  \right]$ for all $k$ follows from the fact that the expectation (i.e. $\mathbf{E}\left[  J_{k}\left(  s_{k}\right)  \right]$) is simply a convex combination of continuous functions. This completes the proof.
\Halmos
\endproof

\smallskip
{\bf Proof of Theorem \ref{2ptsymm}:}
\proof{Proof.}
The results of Theorem \ref{storageTHM} establish that for any finite $N$, the optimal policy is a threshold policy regardless of the price distribution. It can be verified that the form of the optimal policy extends to the infinite-horizon case in the sense that the infinite-horizon optimal policy is similar to the policy defined in Theorem \ref{storageTHM}, with stationary thresholds that are only a function of the storage state. \noindent Furthermore, the optimal average cost is independent of the initial state. Hence, without loss of generality we may assume that we start from $s_0=0$, and therefore, the storage state would only take values that are integer multiples of the ramp rate, i.e., $s\in\{i\overline{v}~|~ i \in 0,1,\cdots,n\}.$ This assumption simplifies the optimal policy to either buy as much as $\overline{v}$, or do nothing, or sell as much as $\overline{v}$. 

\begin{figure}
\[
\xymatrix{
*+[o][F-]{0}\ar@(ul,dl)[]|{a} \ar@/^/[r]^{1-a} & *+[o][F-]{1} \ar@/^/[r]^{1-a} \ar@/^/[l]^{a}& \ldots \ar@/^/[r]^{1-a} \ar@/^/[l]^{a}& *+[o][F-]{ n} \ar@/^/[l]^{a} \ar@(ur,dr)[]|{1-a}\\
}
\]
\caption{Underlying Markov chain associated with storage state under two-point price distribution}
\label{fig:Markov}
\end{figure}
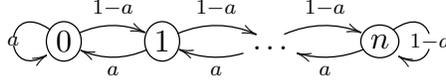

\noindent For conciseness and ease of notation, and without loss of generality, we present the remainder of the proof for the case of $\overline{v}=1, \lambda_{\min}=0$, and $\lambda_{\max}=1$. The proof for the general case is similar. \noindent Let $\overline{t}(s)$ and $\underline{t}(s)$ denote the threshold associated with selling and buying, respectively, at state $s$. Thus, if the stage price $\lambda_k$ falls within the interval $(\underline{t}(s)$, $\overline{t}(s))$, the decision is to do nothing and no cost is incurred. Let $J_k$ denote the cost at stage $k$. Then:
\begin{align}
\mathbf{E}[J_k ]=\mathbf{E}\left[\mathbf{E}[J_k | s_k]\right]&=\mathbf{E}[\mathbf{E}[\lambda_k \min\{\overline{s}-s_k,1\} |s_k, \lambda_k\leq \underline{t}(s_k)]\mathbf{P}(\lambda_k\leq \underline{t}(s_k))\nonumber\\&+\mathbf{E}[\lambda_k\max\{-s_k,-1\} |s_k, \lambda_k\geq\overline{t}(s_k)]\mathbf{P}(\lambda_k\geq\overline{t}(s_k))]\nonumber\\
&\geq \mathbf{E}[\mathbf{E}[\lambda_k\max\{-s_k,-1\} |s_k, \lambda_k\geq\overline{t}(s_k)]\mathbf{P}(\lambda_k\geq\overline{t}(s_k))]]\nonumber\\
&\geq -\mathbf{P}( s_k\geq1)\mathbf{E}[\mathbf{P}(\lambda_k\geq\overline{t}(s_k))]\nonumber\\
&= -(1-P(s_k=0))\mathbf{P}(\lambda_{k}=1) \tag{A.5}\label{LowerB}
\end{align}
\noindent where the first inequality follows from $\lambda_{\min}=0,$ and the second inequality follows from $\lambda_{\max}=1$ and the fact that selling can happen only when $i \neq 0.$ Thus, $\mathbf{E}[J_k ]$ is bounded from below by \eqref{LowerB}, with exact equality holding for a two-point distribution with PMF:
\begin{equation*}
P_{\Lambda}(\lambda)=\left\{ 
\begin{array}
[c]{lcl}
a &\text{if}& \lambda = \lambda_{\max}=1 \\[0.05in]
1-a &\text{if}& \lambda = \lambda_{\min}=0 \\[0.05in]
0 & \text{otherwise}
\end{array}%
\right. 
\end{equation*}
\noindent Next, we show that $a=0.5$ minimizes the stage cost. Under the two-point distribution defined above, the storage state evolves according to the discrete time finite-state Markov chain shown in Figure \ref{fig:Markov}.
%\footnotesize

\noindent Letting $p_i$ denote the steady state probability of state $i,$  
it can be verified that $p_0=1/(1+b+\cdots+b^{n})$, where $b=(1-a)/a$.
\noindent Considering that the optimal decision is to either purchase at a price of $\lambda_{\min}=0$ (which incurs a cost of zero) or sell at a price of $\lambda_{\max}=1$ (which occurs with probability $a(1-p_0)$ and incurs a cost of $-\lambda_{\max}\overline{v}=-1$), we have 
\begin{equation*}
\begin{aligned}
%\min~~
& \gamma =\underset{N\rightarrow \infty}{\lim}
& & \!\!\frac{1}{N}~  \mathbf{E}\left[
{\displaystyle\sum\limits_{k=0}^{N-1}} \lambda_{k}v_k ~|~ s_0=0\right]\\
&~~=\underset{N\rightarrow \infty}{\lim}
& & \!\!\frac{1}{N}~ 
{\displaystyle\sum\limits_{k=0}^{N-1}}\mathbf{E}[ \mathbf{E}\left[ \lambda_{k}\max\{-s_k,-1\} ~|~ \lambda_{k}\geq\overline{t}(s_k), s_k\right]\mathbf{P}(\lambda_k\geq\overline{t}(s_k))]\\
& ~~=-\underset{N\rightarrow \infty}{\lim}
& & \!\!\frac{1}{N}~ 
{\displaystyle\sum\limits_{k=0}^{N-1}} \mathbf{P}(s_k\geq1)\mathbf{P}(\lambda_k=1)=-\underset{N\rightarrow \infty~}{\lim}
 \!\!\frac{1}{N}~  
{\displaystyle\sum\limits_{k=0}^{N-1}}(1-p_0)a,\\
\end{aligned}
%J_\pi(s)=
\label{steadycost}
\end{equation*}
which yields
\begin{equation}
 \gamma =- \frac{b(1+b+\cdots+b^{n-1})}{(b+1)(1+b+\cdots+b^{n})}
 \tag{A.6} \label{gammb}
\end{equation}
\noindent Solving for $b$ to minimize $\gamma$ yields $b=1$, which corresponds to $a=1/2$. This proves \eqref{2ptgamma}, \eqref{2ptdist}, and \eqref{gengamma}. It remains to show that the differential cost function corresponding to \eqref{2ptgamma} is of the form \eqref{diffcost}. 

Let $H(s)$ be the value of the right hand side (RHS) of equation (\ref{valiter}) obtained by plugging the proposed solution for $H^{\ast}(s)$ and $\gamma^{\ast}$ into equation (\ref{valiter}); then we have
\begin{equation}
\begin{aligned}
&H(s)=\mathbf{E} \left[\min_{v \in [\max(-s,-\overline{v}),\min(\overline{v},\overline{s}-s)]}\!\!\!\!\!\!\!\!\!\!\! \lambda v -\frac{(i+1)\lambda_{\min}+(n-i)\lambda_{\max}}{n+1}(s+v)-\frac{i(i+1)(\lambda_{\max}-\lambda_{\min})\overline{v}}{2(n+1)}\right]\\&~~~~~~~~~~ +\frac{\overline{v}(\lambda_{\max}-\lambda_{\min})}{2}\frac{n}{n+1},
\end{aligned}\vspace{-5pt}  \tag{A.7} \label{iter}
\end{equation}
where $i  \overline{v}\leq s+v<(i+1)  \overline{v},~\text{for all} ~i\in\{0,1,2,\cdots,n-1\}$. 

\noindent Therefore, we need to show that $H(s)$, obtained from solving the optimization problem in (\ref{iter}), is indeed the $H^{\ast}(s)$ defined in (\ref{diffcost}). It can be verified that the solution to the optimization problem in (\ref{iter}) is as follows:
%$~
\begin{equation*}
v^{\ast}=\left\{
\begin{array}
[c]{lcl}%
\max(-s,-\overline{v}) & \text{if} & \lambda=\lambda_{\max}\\[0.06in]
\min(\overline{v},\overline{s}-s) & \text{if} & \lambda=\lambda_{\min} \\
\end{array}
\right.
 \tag{A.8} \label{bellopt}
\end{equation*}

\noindent It can be verified by inspection that if we plug (\ref{bellopt}) into (\ref{iter}), the $H(s)$ obtained is indeed the $H^{\ast}(s)$ defined in (\ref{diffcost}). The proof is complete.
\Halmos
\endproof
\smallskip
{\bf Proof of Corollary \ref{meancoroll}:}
\proof{Proof.}
Since a two-point distribution with non-zero probability masses placed at the endpoints of the fixed support is always achievable for any $\mu \in (\lambda_{\min},\lambda_{\max})$, the additional assumption of a fixed mean does not affect the proof of Theorem \ref{2ptsymm} up to and including equation \eqref{gammb}. For the proof of this corollary, we need to minimize $\gamma$ in equation \eqref{gammb} with respect to $b$ subject to the constraint $\mu= b\lambda_{\min}/(b+1)+\lambda_{\max}/(b+1)$. This optimization problem has the unique solution $b=(\lambda_{\max}-\mu)/(\mu-\lambda_{\min})$. Recall that $a=1/(b+1)$. Hence, $a=(\mu-\lambda_{\min})/(\lambda_{\max}-\lambda_{\min})$. This completes the proof.
\Halmos
\endproof

\end{document}